\documentclass[11pt]{amsart}

\usepackage{amsmath}
\usepackage{amssymb}
\usepackage{amsthm}
\usepackage{enumitem}
\usepackage{microtype}
\usepackage{array}
\usepackage{dcolumn}
\usepackage[margin=1.0in]{geometry}
\usepackage{float}
\usepackage{graphicx}
\usepackage{tikz-cd}
\usepackage[style=alphabetic]{biblatex}
\usepackage[hidelinks]{hyperref}
\usepackage{mathtools}
\usepackage{todonotes}
\usepackage{mathrsfs}

\renewcommand \v{\vert}

\renewcommand \Im{\text{Im}}

\newcommand \R{\mathbb{R}}

\newcommand \Z{\mathbb{Z}}

\newcommand \ep{\varepsilon}
\newcommand \la{\langle}
\newcommand \ra{\rangle}
\newcommand \del{\partial}

\DeclareMathOperator \Ker{Ker}
\DeclareMathOperator \Int{Int}

\DeclareMathOperator \RP{\mathbb{RP}}

\DeclareMathOperator \Aut{Aut}
\DeclareMathOperator \Mod{Mod}
\DeclareMathOperator \PMod{PMod}
\DeclareMathOperator \GL{GL}

\DeclareMathOperator \SO{SO}

\DeclareMathOperator \id{id}

\DeclareMathOperator \LMod{LMod}
\DeclareMathOperator \SMod{SMod}
\DeclareMathOperator \Twist{Twist}

\DeclareMathOperator \Diff{Diff}
\DeclareMathOperator \Out{Out}
\DeclareMathOperator \Inn{Inn}
\DeclareMathOperator \SymOut{SymOut}
\DeclareMathOperator \PSymOut{PSymOut}
\DeclareMathOperator \SymAut{SymAut}
\DeclareMathOperator \PSymAut{PSymAut}
\DeclareMathOperator \PalAut{PalAut}

\DeclareMathOperator \st{st}
\DeclareMathOperator \Isom{Isom}
\DeclareMathOperator \Homeo{Homeo}

\DeclareMathOperator \rel{rel}
\DeclareMathOperator \Fr{Fr}
\DeclareMathOperator \UT{UT}
\DeclareMathOperator \Emb{Emb}
\renewcommand \L{\mathcal{L}}

\newcommand \F{\mathcal{F}}

\renewcommand \P{\mathcal{P}}
\newcommand \T{\mathscr{T}}
\newcommand \B{\mathcal{B}}
\newcommand \calS{\mathcal{S}}
\newcommand \K{\mathcal{K}}

\theoremstyle{definition}
\newtheorem{definition}{Definition}[section]

\theoremstyle{plain}
\newtheorem{question}[definition]{Question}

\newtheorem{lemma}[definition]{Lemma}
\newtheorem{proposition}[definition]{Proposition}

\newtheorem{mainTheorem}{Theorem}

\newtheorem{mainCorollary}[mainTheorem]{Corollary}
\theoremstyle{remark}
\newtheorem{remark}[definition]{Remark}

\numberwithin{equation}{section}

\title{Birman-Hilden theory for 3-manifolds}
\author{Trent Lucas}

\begin{document}

\begin{abstract}
    Given a branched cover of manifolds, one can lift homeomorphisms along the cover to obtain a (virtual) homomorphism between mapping class groups.
    Following a question of Margalit-Winarski, we study the injectivity of this lifting map in the case of $3$-manifolds.
    We show that in contrast to the case of surfaces, the lifting map is generally not injective for most regular branched covers of $3$-manifolds.
    This includes the double cover of $S^3$ branched over the unlink, which generalizes the hyperelliptic branched cover of $S^2$.
    In this case, we find a finite normal generating set for the kernel of the lifting map.
\end{abstract}

\maketitle

\section{Introduction}\label{subsec:intro}
\subsection{Motivating questions and main results}
Given a branched cover of manifolds, one obtains an associated \emph{lifting map} between mapping class groups.  
Namely, suppose $p:M \rightarrow N$ is a finite branched cover of closed oriented manifolds with deck group $G$.
Let $C \subseteq N$ be the branch set of $p$, and let $\Mod(N,C)$ be the mapping class group of the pair $(N,C)$.
This is the group $\pi_0(\Homeo^+(N,C))$, where $\Homeo^+(N,C)$ is the group homeomorphisms of $N$ that preserve $C$ setwise and preserve the orientation of $N$ (but may reverse the orientation of the components of $C$).
Let $\Gamma$ denote the image of $G$ in the mapping class group $\Mod(M) \coloneqq \pi_0(\Homeo^+(M))$.
Then sending a mapping class to its lift along $p$ yields a {lifting map}
\begin{equation*}
    \L_p:\LMod_p(N) \rightarrow \SMod_p(M)/\Gamma.
\end{equation*}
Here $\LMod_p(N) \leq \Mod(N, C)$ is the finite-index subgroup of liftable mapping classes, and $\SMod_p(M) \leq \Mod(M)$ is the subgroup of mapping classes represented by a \emph{symmetric} homeomorphism, i.e.\ a homeomorphism $f$ such that $p(x) = p(y)$ implies $p(f(x)) = p(f(y))$.

This paper is motivated by the following question, which was posed (in a slightly different form) by Margalit-Winarski \cite[Ques~11.4]{margalit-winarski}.

\begin{question}\label{ques:injectivity}
    For which branched covers of $3$-manifolds is the lifting map $\L_p$ injective?
\end{question}

A theorem of Birman-Hilden \cite{birman-hilden} and MacLachlan-Harvey \cite{maclachlan-harvey} says that when $M$ is a hyperbolic surface and $p$ is regular, $\L_p$ is injective.
This result has had numerous applications for both surface mapping class groups and $3$-manifold topology; see Margalit-Winarski's survey \cite{margalit-winarski} for examples.
In contrast to the surface case, we show that the lifting map $\L_p$ is \emph{not} injective for most finite regular branched covers of $3$-manifolds. 

\begin{mainTheorem}\label{mainthm:BH_counterexamples}
    Let $p:M \rightarrow N$ be a finite regular branched cover of closed oriented $3$-manifolds given by a finite group action $G \leq \Diff^+(M)$.
    Assume the branch set $C \subseteq N$ is nonempty, and let $N^\circ$ denote the complement of a regular open neighborhood of $C$.
    If either
    \begin{enumerate}[label=(\roman*)]
        \item $N^\circ$ has at least $3$ prime factors, or
        \item $N^\circ$ has $2$ prime factors, and at least one prime factor $W$ has nonempty boundary and is not a Seifert fibration over the disk,
    \end{enumerate}
    then the kernel of the lifting map $\L_p:\LMod_p(N) \rightarrow \SMod_p(M)/\Gamma$ is infinite.
\end{mainTheorem}

Theorem \ref{mainthm:BH_counterexamples} implies that most regular branched covers of $3$-manifolds do not have the \emph{Birman-Hilden property} as defined by Margalit-Winarski \cite{margalit-winarski}.
This means that for any cover $p:M \rightarrow N$ as in Theorem \ref{mainthm:BH_counterexamples}, there is a pair of symmetric homeomorphisms of $M$ that are isotopic, but not isotopic through symmetric homeomorphisms.
See Remark \ref{rem:BH_prop_vs_injectivity} for the precise relationship between $\L_p$ and the Birman-Hilden property.

Theorem \ref{mainthm:BH_counterexamples} is sharp in the following sense.
The cases where $N^\circ$ is reducible that are not covered by Theorem \ref{mainthm:BH_counterexamples} are when $N^\circ \cong W_1^\circ \# W_2^\circ$, where $W_1^\circ$ is a Seifert fibration over the disk and $W_2^\circ$ is another such fibration or a closed prime $3$-manifold.
In these cases, $\Ker(\L_p)$ may not be infinite.
The simplest example is the double cover $p_2:S^1 \times S^2 \rightarrow S^3$ defined below, where each $W_i^\circ$ is a solid torus.
In this example, $\LMod_{p_2}(S^3)$ is finite, and hence $\Ker(\L_{p_2})$ is too (we describe $\Ker(\L_{p_2})$ explicitly below).
On the other hand, there are examples with $W_1^\circ$ a solid torus and $W_2^\circ$ closed where $\Ker(\L_p)$ is infinite (see Section \ref{subsec:exceptional_cases}).
It appears to be a subtle problem to determine when $\Ker(\L_p)$ is infinite in the reducible cases not covered by Theorem \ref{mainthm:BH_counterexamples}.
For a discussion of the case that $N^\circ$ is irreducible, see Section \ref{subsec:irreducible_case} below.

Theorem \ref{mainthm:BH_counterexamples} puts the hypotheses on the topology of the manifold $N^\circ$.
One might instead hope to understand the injectivity of $\L_p$ in terms of the action of $G$ on $M$.
From part (i) of Theorem \ref{mainthm:BH_counterexamples}, we deduce such a statement using the Equivariant Sphere Theorem of Meeks-Yau \cite{meeks-yau-sphere}.

\begin{mainCorollary}\label{maincor:BH_counterexamples}
    Assume $G$ does not act freely on $M$, and let $M^\circ$ denote the complement of a regular open neighborhood of $p^{-1}(C)$.
    If the group of coinvariants $\pi_2(M^\circ)_G$ has rank at least $3$ as a $\pi_1(M^\circ)$-module, then $\Ker(\L_p)$ is infinite.
\end{mainCorollary}

Given that $\L_p$ is frequently non-injective for $3$-manifolds, the natural next goal is to understand its kernel.
Perhaps the first question one would hope to answer is the following.

\begin{question}\label{ques:kernel}
    If the lifting map $\L_p$ is not injective, what is a natural generating set for its kernel?
\end{question}

We expect that Question \ref{ques:kernel} is difficult to answer in general.
The kernel of $\L_p$ is generally an infinite index subgroup of $\Mod(N,C)$, and hence need not share any finiteness properties that $\Mod(N,C)$ may have.
Regardless, we answer Question \ref{ques:kernel} for an infinite family of covers.

Namely, let $n \geq 2$, and let $p_n:M_n \rightarrow S^3$ be the double cover branched over the $n$-component unlink $C_n \subseteq S^3$.
The manifold $M_n$ is the $(n-1)$-fold connected sum $(S^1 \times S^2)^{\#(n-1)}$, and the deck group $G$ is generated by an order $2$ diffeomorphism $\tau$ that reverses the orientation of each $S^1$ and $S^2$ factor.
In this case, $\LMod_{p_n}(S^3)$ is the full mapping class group $\Mod(S^3, C_n)$; this group is closely related to the ``loop braid group'' (see e.g.\ \cite{damiani}).
We denote the lifting map by
\begin{equation*}
    \L_n:\Mod(S^3, C_n) \rightarrow \SMod_n(M_n)/\la \tau \ra.
\end{equation*}
The cover $p_n$ is a direct analog of the hyperelliptic involution on a genus $n-1$ surface.

Margalit-Winarski discuss the cover $p_n$ specifically, and ask whether $\L_n$ is injective \cite[Ques~11.5]{margalit-winarski}.
Our work shows that the answer is ``no.''
In the case $n = 2$, the groups $\Mod(S^3,C_2)$ and $\Mod(S^1 \times S^2)$ are finite, and we describe $\Ker(\L_2)$ explicitly below (see also Section \ref{subsec:two_factor_case}).
For $n \geq 3$, Theorem \ref{mainthm:BH_counterexamples} implies that $\Ker(\L_n)$ is infinite.
Our second result gives a normal generating set in these cases.

\begin{mainTheorem}\label{mainthm:hyperelliptic_kernel}
    Let $\rho \in \Mod(S^3, C_n)$ be the order $2$ element which reverses the orientation of each component of $C_n$.
    Then for $n \geq 3$, $\Ker(\L_n)$ is the $\Mod(S^3,C_n)$-normal closure of $\{\rho\}$.
\end{mainTheorem}

It is not hard to show that $\rho$ lies in $\Ker(\L_n)$.
Indeed, $\rho$ lifts to a product of sphere twists which is well-known to be trivial; see Lemma \ref{lem:kernel_lower_bound} for details.
Thus, Theorem \ref{mainthm:hyperelliptic_kernel} says that $\Ker(\L_n)$ is as small as possible given the ``obvious'' element in the kernel.

A theorem of Goldsmith \cite{goldsmith}, building off the work of Dahm \cite{dahm}, identifies $\Mod(S^3, C_n)$ as the group of symmetric outer automorphisms of $F_n$ (see Section \ref{subsec:alg_perspective} below).
This group has a simple finite presentation by the work of McCool \cite{mccool}, as well as Gilbert \cite{gilbert} following Fouxe-Rabinovitch \cite{fouxe-rabinovitch}.
Theorem \ref{mainthm:hyperelliptic_kernel} therefore allows one to obtain a presentation for $\SMod_n(M_n)$.
We can also use the presentation of $\Mod(S^3, C_n)$ and Theorem \ref{mainthm:hyperelliptic_kernel} to explicitly identify $\Ker(\L_n)$ for small $n$.

\begin{mainCorollary}\label{maincor:explicit_kernel_computation}
    For $n \in \{2,3\}$, $\Ker(\L_n)$ is as follows.
    \begin{enumerate}[label=(\roman*)]
        \item For $n=2$, $\Ker(\L_2) \cong \Z/2\Z \times \Z/2\Z$, and is generated by $\rho$ and the element $\sigma$ which swaps the two components of $C_2$.
        \item For $n=3$, $\Ker(\L_3)$ splits as a semidirect product $F_\infty \rtimes \Z/2\Z$, where $F_\infty$ is an infinite rank free group.  In particular, $\Ker(\L_3)$ is not finitely generated.
    \end{enumerate}
\end{mainCorollary}

Corollary \ref{maincor:explicit_kernel_computation} is reminiscent of the Torelli group $\mathcal{I}_g$, which is trivial for genus $g=1$ and an infinite rank free group for $g=2$ \cite{mess}.
It would be of interest to determine whether $\Ker(\L_n)$ is finitely generated for $n \geq 3$, as $\mathcal{I}_g$ is for $g \geq 3$ \cite{johnson}.

\subsection{Application: trivial bundles with a nontrivial fiberwise action}
Theorem \ref{mainthm:BH_counterexamples} provides many examples of a $3$-manifold bundle $E \rightarrow S^1$ equipped with a fiberwise group action such that $E$ is topologically trivial, but not equivariantly trivial.
Let $M$ be a closed oriented $3$-manifold and $G \leq \Diff^+(M)$ a finite group, and suppose that the quotient map $M \rightarrow M/G$ satifies the hypotheses of Theorem \ref{mainthm:BH_counterexamples}.
Let $\Homeo_G^+(M) \leq \Homeo^+(M)$ denote the subgroup of homeomorphisms that commute with $G$.
Then Theorem \ref{mainthm:BH_counterexamples} implies that the map $\pi_0(\Homeo_G^+(M)) \rightarrow \pi_0(\Homeo^+(M))$ is not injective (see Remark \ref{rem:BH_prop_vs_injectivity}).
This implies that the natural map 
\begin{equation*}
    \kappa:\pi_1(B\Homeo^+_G(M)) \rightarrow \pi_1(B\Homeo^+(M))
\end{equation*}
is not injective, where $B$ denotes the classifying space.
Elements of $\pi_1(B\Homeo^+(M))$ correspond to isomorphism classes of $M$-bundles over $S^1$.
Elements of $\pi_1(B\Homeo^+_G(M))$ correspond to isomorphism classes of $M$-bundles over $S^1$ with structure group $\Homeo^+_G(M)$; these are bundles for which the action of $G$ on $M$ extends to a fiberwise action on the total space.
Thus, elements of $\Ker(\kappa)$ correspond to bundles which are topologically trivial, but not equivariantly trivial, i.e.\ not equivariantly isomorphic to the product $S^1 \times M$ equipped with the standard action of $G$ in each fiber. 

This application suggests the following question, which is a ``higher'' version of Question \ref{ques:injectivity}.

\begin{question}\label{ques:higher_BH}
    Suppose a finite group $G$ acts on a closed oriented $3$-manifold $M$.
    For $k \geq 1$, is the map $\pi_k(\Homeo_G^+(M)) \rightarrow \pi_k(\Homeo^+(M))$ injective?
\end{question}

As above, if the answer is ``no,'' then we obtain $M$-bundles over $S^{k+1}$ which are topologically but not equivariantly trivial.

\subsection{The algebraic perspective of Theorem \ref{mainthm:hyperelliptic_kernel}}\label{subsec:alg_perspective}
Recall that $\pi_1(S^3 - C_n)$ is a free group on $n$ generators.
The action of $\Mod(S^3,C_n)$ on $\pi_1(S^3 - C_n)$ therefore induces a map $\Mod(S^3,C_n) \rightarrow \Out(F_n)$.
As mentioned above, a theorem of Goldsmith \cite{goldsmith} says that this map yields an isomorphism 
\begin{equation*}
    \Mod(S^3,C_n) \cong \SymOut(F_n).
\end{equation*}
Here $\SymOut(F_n)$ is the group of \emph{symmetric} outer automorphisms; this is the group of outer automorphisms that map each generator to a conjugate of another generator or its inverse.

Moreover, recall that $\pi_1(M_n) \cong F_{n-1}$. 
A theorem of Brendle-Broaddus-Putman \cite{brendle-broaddus-putman}, building off the work of Laudenbach \cite{laudenbach-1,laudenbach-2}, says that the action on $\pi_1$ yields an isomorphism
\begin{equation*}
    \Mod(M_n) \cong (\Z/2\Z)^{n-1} \rtimes \Out(F_{n-1}).
\end{equation*}
Here $(\Z/2\Z)^{n-1}$ is the subgroup of sphere twists, which act trivially on $\pi_1$.
The involution $\tau$ induces the outer automorphism of $F_{n-1}$ that inverts each generator.
Consequently, the group $\SMod_n(M_n)$ is closely related to the group of \emph{palindromic} automorphisms of $F_{n-1}$; see Remark \ref{rem:palindromic_automorphisms} for details.

The lifting map $\L_n:\Mod(S^3,C_n) \rightarrow \SMod_n(M_n)/\la \tau \ra$ therefore yields a virtual homomorphism $\SymOut(F_n) \rightarrow \Out(F_{n-1})$. 
It is not a priori obvious how one would define such a map in algebraic terms.
However, the algebraic perspective is clarified by the \emph{relative fundamental group} defined below; using the relative fundamental group, we show that the lifting map $\L_n$ corresponds to the mod $2$ reduction map $\SymOut(F_n) \rightarrow \SymOut((\Z/2\Z)^{*n})$ (see Section \ref{subsec:pf_idea_thm_C}).

\subsection{Key idea: the relative fundamental group}

Our main approach in this paper is to understand mapping classes by their action on $\pi_1$.
For branched covers, this can be tricky.
Unlike the unbranched case, if the cover $p:M \rightarrow N$ is branched, then $\pi_1(M)$ is not naturally a subgroup of $\pi_1(N)$.
Instead, if we let $M^\circ$ and $N^\circ$ denote the branch set complements, then $\pi_1(M)$ is a quotient of $\pi_1(M^\circ)$, which is a subgroup of $\pi_1(N^\circ)$.
This makes it hard to read off the action of $\LMod_p(N)$ on $\pi_1(M)$ directly from its action on $\pi_1(N)$ or $\pi_1(N^\circ)$.

Instead, we work with the \emph{relative fundamental group}.
Fix a base point $* \in M - p^{-1}(C)$.
We define $\pi_1^{\rel}(M)$ to be the relative homotopy ``group'' $\pi_1(M, G*, *)$.
Using the action of $G$, $\pi_1^{\rel}(M)$ has a natural group structure (see Section \ref{subsec:group_structure}), and it fits into the following commutative diagram with exact rows:
\begin{equation*}
    \begin{tikzcd}
        1 & {\pi_1(M^\circ)} & {\pi_1(N^\circ)} & G & 1 \\
        1 & {\pi_1(M)} & {\pi_1^{\rel}(M)} & {\pi_0(G*)} & 1
        \arrow[from=1-1, to=1-2]
        \arrow[from=1-2, to=1-3]
        \arrow[two heads, from=1-2, to=2-2]
        \arrow[from=1-3, to=1-4]
        \arrow[two heads, from=1-3, to=2-3]
        \arrow[from=1-4, to=1-5]
        \arrow["\cong", from=1-4, to=2-4]
        \arrow[from=2-1, to=2-2]
        \arrow[from=2-2, to=2-3]
        \arrow[from=2-3, to=2-4]
        \arrow[from=2-4, to=2-5]
    \end{tikzcd}
\end{equation*}
Here the top row comes from the monodromy homomorphism of the cover $M^\circ \rightarrow N^\circ$, and the bottom row comes from the long exact sequence of the pair $(M, G*)$.

The upshot is that we obtain the following commutative diagram (see Proposition \ref{prop:lifting_map_factors}):
\begin{equation*}
    \begin{tikzcd}
        {\LMod_p(N)} & {\Out(\pi_1^{\rel}(M))} \\
        {\SMod_p(M)/\Gamma} & {\Out_{\Gamma_*}(\pi_1(M))/\Gamma_*}
        \arrow["\Psi", from=1-1, to=1-2]
        \arrow["{\L_p}"', from=1-1, to=2-1]
        \arrow["\overline{r}", dashed, from=1-2, to=2-2]
        \arrow["\Theta", from=2-1, to=2-2]
    \end{tikzcd}
\end{equation*}
Here $\Gamma_*$ is the image of $G$ in $\Out(\pi_1(M))$, $\Out_{\Gamma_*}(\pi_1(M))$ is its normalizer, and $\Theta$ is the natural map.
The maps $\Psi$ and $\overline{r}$ are defined in Section \ref{sec:rel_fundamental_group} (the dashed arrow $\overline{r}$ is defined only on the image of $\Psi$).
The kernel of the map $\Theta$ is finite (see Proposition \ref{prop:lifting_map_factors}).
Therefore, this diagram lets us find elements of $\Ker(\L_p)$ by instead finding infinite order elements of $\Ker(\Psi)$, which proves to be easier.
Moreover, if the map $\overline{r}$ is injective, then computing $\Ker(\L_p)$ essentially amounts to computing $\Ker(\Psi)$; this occurs for the cover $p_n:M_n \rightarrow S^3$, where the map $\Psi$ has a simple algebraic description as we explain below.

The relative fundamental group is closely related to the orbifold fundamental group (see Remark \ref{rem:orbifold_pi1}).
While the relative fundamental group is defined only for regular covers, its advantage over the orbifold fundamental group is that it is simpler to define and work with for the purpose of proving Theorems \ref{mainthm:BH_counterexamples} and \ref{mainthm:hyperelliptic_kernel}.

\subsection{Proof idea of Theorem \ref{mainthm:BH_counterexamples}}
To prove Theorem \ref{mainthm:BH_counterexamples}, we find an infinite order element in the kernel of the map $\Psi:\LMod_p(N) \rightarrow \Out(\pi_1^{\rel}(M))$.
The element we find is a certain \emph{slide homeomorphism}, which acts on $\pi_1(N^\circ)$ by \emph{partial conjugation}.
That is, a connected sum decomposition $N^\circ \cong W_1^\circ \# W_2^\circ$ induces a free splitting $\pi_1(N^\circ) \cong \pi_1(W_1^\circ) * \pi_1(W_2^\circ)$; a slide homeomorphism will conjugate one free factor and fix the other pointwise.
The free splitting of $\pi_1(N^\circ)$ does not necessarily induce a free splitting of $\pi_1(M)$, so it can be difficult to determine how a lift of a slide homeomorphism acts on $\pi_1(M)$.
However, the key insight is that the free splitting of $\pi_1(N^\circ)$ does descend to a free splitting $\pi_1^{\rel}(M) \cong H_1 * H_2$ (see Proposition \ref{prop:rel_pi1_is_free_product}), and a lift of a slide homeomorphism will act by partial conjugation on $\pi_1^{\rel}(M)$.
Thus, to find an infinite order element of $\Ker(\Psi)$, we must find an element $\gamma \in \pi_1(W_i^\circ)$ that has no central power, but lifts to an element which is central in $H_i$.

For case (i), we show that one can take $\gamma$ to be a meridian around the branch set in $N^\circ$.
Case (ii) is more difficult, because such a meridian may have a central power.
If this happens, the Seifert Fiber Space Theorem \cite{waldhausen-SFST} implies that the summand $W_i^\circ$ containing this meridian is Seifert fibered.
Thus, to prove Theorem \ref{mainthm:BH_counterexamples} in case (ii), we construct elements in the lifting kernel for certain branched covers of Seifert fibered manifolds (see Lemma \ref{lem:SF_kernel_infinite}).

Slide homeomorphisms can also be defined for surfaces, but they behave quite differently in the $3$-manifold case compared to the surface case.
This difference is rooted in the fact that sphere twists on $3$-manifolds behave much differently from Dehn twists on surfaces (namely, sphere twists have order $2$ and act trivially on $\pi_1$).
This difference helps explain why slide homeomorphisms provide elements in the kernel of the lifting map for $3$-manifolds, but not for surfaces.
See Remark \ref{rem:surface-slides} for more details.

\subsection{Proof idea of Theorem \ref{mainthm:hyperelliptic_kernel}}\label{subsec:pf_idea_thm_C}
Recall from above that for the double cover $p_n:M_n \rightarrow S^3$ branched over the unlink $C_n \subseteq S^3$, the mapping class group $\Mod(S^3, C_n)$ is isomorphic to the group $\SymOut(F_n)$ of symmetric outer automorphisms of $F_n$.
We show that in this case, $\pi_1^{\rel}(M_n)$ is isomorphic to the $n$-fold free product $(\Z/2\Z)^{*n}$, which we denote by $H_n$ (see Lemma \ref{lem:hyperelliptic_rel_htpy_groups}).
Thus, the map $\Psi:\Mod(S^3, C_n) \rightarrow \Out(\pi_1^{\rel}(M))$ is the natural map
\begin{equation*}
    \P_n:\SymOut(F_n) \rightarrow \SymOut(H_n)
\end{equation*}
obtained by reducing each generator of $F_n$ mod $2$.
We then reduce the problem of computing $\Ker(\L_n)$ to the purely algebraic problem of computing $\Ker(\P_n)$.
We show the following, which may be of independent interest (see Section \ref{sec:alg_kernel}).

\begin{proposition}
    For $1 \leq i \leq n$, let $\rho_i \in \SymOut(F_n)$ be the element represented by the automorphism of $F_n$ that inverts the $i$th generator.
    Then the kernel of the map $\P_n:\SymOut(F_n) \rightarrow \SymOut(H_n)$ is the $\SymOut(F_n)$-normal closure of the set $\{\rho_1, \ldots, \rho_n\}$.
\end{proposition}

We prove this proposition using the action of symmetric outer automorphism groups on certain contractible simplicial complexes.
For any free product of indecomposable groups $L = L_1 * \cdots * L_n$, McCullough-Miller \cite{mccullough-miller} constructed a contractible simplicial complex $K_0(L)$ on which $\SymOut(L)$ acts.
The key insight is that $\Ker(\P_n)$ acts on the contractible complex $K_0(F_n)$, and the quotient is the contractible complex $K_0(H_n)$ (see Lemma \ref{lem:quotient_of_complexes}).
It follows that $\Ker(\P_n)$ is generated by its vertex stabilizers on $K_0(F_n)$.
We compute these vertex stabilizers using McCullough-Miller's computation of the $\SymOut(F_n)$ vertex stabilizers.

In principle, one could use a result of Putman \cite{putman-choices} to compute an (infinite) presentation of $\Ker(\P_n)$ from the isomorphism $K_0(F_n)/\Ker(\P_n) \cong K_0(H_n)$.
It would be of interest to study this presentation, and perhaps determine whether $\Ker(\P_n)$ is finitely generated or presented for $n \geq 3$.

\subsection{The relative fundamental group and braid groups}\label{subsec:rel_pi1_and_braid}
The relative fundamental group can shed new light on the lifting map for surfaces; here is one example (cf.\ \cite[\S 7]{margalit-winarski}).
Let $D_n$ denote the $n$-times punctured disk, and recall that $\pi_1(D_n)$ is the free group $F_n$.
The action of the braid group $B_n$ on $\pi_1(D_n)$ yields an embedding $B_n \hookrightarrow \SymAut(F_n)$, where $\SymAut(F_n)$ is the group of symmetric automorphisms.
Let $H_{n,k}$ denote the $n$-fold free product $(\Z/k\Z)^{*n}$.
Then there is a natural map $\SymAut(F_n) \rightarrow \SymAut(H_{n,k})$, yielding a map 
\begin{equation*}
    \eta_{n,k}:B_n \rightarrow \SymAut(H_{n,k}).
\end{equation*}
Magnus asked whether the map $\eta_{n,k}$ is injective.
Birman-Hilden \cite{birman-hilden} proved that it is injective using the injectivity of the lifting map for a certain branched cover of the disk $p_{n,k}:S_{n,k} \rightarrow D$ (Magnus later applied this result in \cite{magnus}).

At first glance, Birman-Hilden's proof is a surprising application of the lifting map in a purely algebraic problem.
However, its relevance is explained by the following observation: for the cover $p_{n,k}$, the group $H_{n,k}$ is precisely the relative fundamental group $\pi_1^{\rel}(S_{n,k})$.
Furthermore, in the case $k=2$, a purely algebraic proof that $\eta_{n,k}$ is injective (as in \cite{bacardit-dicks} or \cite{johnson-problem}) yields a proof that the lifting map $B_n \rightarrow \SMod(S_{n,2})$ is injective.
See Remark \ref{rem:rel_pi1_and_braids_discussion} for more details.

\subsection{The irreducible case}\label{subsec:irreducible_case}
There are many natural questions left to answer about the lifting map for $3$-manifolds.
The first is whether one can prove a complementary result to Theorem \ref{mainthm:BH_counterexamples} for the case of irreducible $3$-manifolds.

\begin{question}\label{ques:reducible_case}
    Let $p:M \rightarrow N$ be a finite regular branched cover of $3$-manifolds.
    Suppose that the branch set complement $N^\circ$ is irreducible.
    Is the lifting map $\L_p$ injective?
\end{question}

The answer to Question \ref{ques:reducible_case} is ``no'' in some cases.
Namely, as part of the proof of Theorem \ref{mainthm:BH_counterexamples}, we show in Lemma \ref{lem:SF_kernel_infinite} that $\Ker(\L_p)$ is infinite for certain covers where $N^\circ$ is Seifert fibered.
It would be especially interesting if one could answer Question \ref{ques:reducible_case} for Seifert fibrations over the disk, since then one might be able to understand the cases not covered by Theorem \ref{mainthm:BH_counterexamples}.
We also note that Birman-Hilden and MacLachlan-Harvey's proofs of the injectivity of $\L_p$ in the surface case use hyperbolic geometry in a crucial way.
It is therefore natural to ask in particular about the case that $N^\circ$ is hyperbolic.

\subsection{The unbranched case}
All of our results concern branched covers.
It is of course natural to ask about the potentially simpler case of unbranched covers.

\begin{question}\label{ques:unbranched_case}
    Let $p:M \rightarrow N$ be a finite unbranched cover of $3$-manifolds.
    Is the lifting map $\L_p$ injective?
\end{question}

The answer to Question \ref{ques:unbranched_case} is ``yes'' in certain cases.
Vogt \cite{vogt} proved that $\L_p$ is injective for certain unbranched covers where $N$ is Seifert fibered.
Moreover, if $N$ is aspherical and roots are unique in $\pi_1(N)$ (such as when $N$ is hyperbolic or $N = T^3$), and if homotopic homeomorphisms of $N$ are always isotopic (such as when $N$ is irreducible and Haken \cite[Thm~7.1]{waldhauden-pi1-iso}), then $\L_p$ is injective.
This follows directly from an argument in the surface case due to Birman-Hilden \cite{birman-hilden-unbranched-argument} and Aramayona-Leininger-Souto \cite{aramayona-leininger-souto}; we follow \cite[\S 9]{margalit-winarski}.
The argument is as follows: if $\alpha \in \Ker(\L_p)$, then $\alpha$ has a representative $f$ that acts on $\pi_1(N)$ and pointwise fixes the subgroup $\pi_1(M) \leq \pi_1(N)$.
Since roots are unique, this implies that $f$ acts trivially on $\pi_1(N)$.
Since $f$ is aspherical, this implies that $f$ is homotopic to the identity, and thus $f$ is isotopic to the identity.

\subsection{Outline}
In Section \ref{sec:preliminaries}, we recall some basic facts about regular branched covers of $3$-manifolds and formally define the lifting map $\L_p$.
In Section \ref{sec:rel_fundamental_group}, we introduce the relative fundamental group $\pi_1^{\rel}(M)$, and we prove Proposition \ref{prop:lifting_map_factors}, which relates $\Ker(\L_p)$ to the action of $\LMod_p(N)$ on $\pi_1^{\rel}(M)$.
In Section \ref{sec:pf_thm_A}, we prove Theorem \ref{mainthm:BH_counterexamples} and Corollary \ref{maincor:BH_counterexamples}, and discuss the exceptional cases of Theorem \ref{mainthm:BH_counterexamples}.
In Section \ref{sec:he_cover}, we define the cover $p_n:M_n \rightarrow S^3$ and discuss the liftable and symmetric mapping class groups for this cover.
In Section \ref{sec:alg_characterization}, we prove Proposition \ref{prop:lifting_kernel_is_alg_kernel}, which relates $\Ker(\L_p)$ to the kernel of the algebraically defined map $\P_n:\SymOut(F_n) \rightarrow \SymOut(H_n)$.
In Section \ref{sec:alg_kernel}, we compute $\Ker(\P_n)$ using the action of symmetric outer automorphism groups on McCullough-Miller's complexes.
Finally, in Section \ref{sec:pf_thm_C}, we complete the proof of Theorem \ref{mainthm:hyperelliptic_kernel} and Corollary \ref{maincor:explicit_kernel_computation}.

\subsection{Acknowledgements}
We thank Bena Tshishiku for introducing us to this problem, for helpful comments on an earlier draft of this paper, and for many helpful conversations throughout this project.
We thank Tyrone Ghaswala for helpful comments, and in particular for pointing out the connection between the relative fundamental group and the orbifold fundamental group.
We also thank Gabriel Corrigan and Ryan Dickmann for helpful conversations about the concepts in this paper, and we thank Dan Margalit for comments on an earlier draft.

\section{Preliminaries}\label{sec:preliminaries}

Throughout this paper, we fix an action of a finite group $G \leq \Diff^+(M)$ on a closed oriented smooth 3-manifold $M$.
We let $N = M/G$ let $p:M \rightarrow N$ denote the quotient map.
We let $B \subseteq M$ denote the \emph{ramification set}, i.e.\ the set of points with nontrivial $G$-stabilizers, and we let $C \subseteq N$ denote the \emph{branch set}, i.e.\ the image of $B$ under $p$.
We assume throughout this paper that $B$ and $C$ are nonempty unless specified otherwise.
We show below that $B$ and $C$ are embedded graphs with no contractible components, and that there are handlebody neighborhoods $V \subseteq M$ and $\overline{V} \coloneqq p(V) \subseteq N$ which deformation retract onto $B$ and $C$ respectively.
We let $M^\circ$ and $N^\circ$ denote the compact manifolds $M - \Int(V)$ and $N - \Int(\overline{V})$ respectively, so $M^\circ$ and $N^\circ$ have nonempty boundary with no spherical components.
Note that $p$ restricts to an unbranched cover $M^\circ \rightarrow N^\circ$.

In this section, we establish some basic results about the branched cover $p$ and the lifting map $\L_p$.
In Section \ref{subsec:branch_set}, we show that $B$ and $C$ are embedded graphs that admit compatible handlebody neighborhoods. 
In Section \ref{subsec:group_actions_on_handlebodies}, we describe the local behavior of $p$ on these neighborhoods.
In Section \ref{subsec:liftable_and_symmetric_mcgs}, we formally define the liftable and symmetric mapping class groups, and in Section \ref{subsec:lifting_map}, we define the lifting map $\L_p$.

\subsection{The branch set}\label{subsec:branch_set}
In the topological category, the fixed set of a finite-order homeomorphism may be a ``wild'' submanifold.
By contrast, in the smooth category, the fixed set is quite regular.
The following lemma is routine, but we include it for completeness.

\begin{lemma}\label{lem:branch_set_regularity}
    The ramification set $B \subseteq M$ is a (possibly disconnected) embedded finite graph with no contractible components, and there is a $G$-invariant closed handlebody $V \subseteq M$ that deformation retracts onto $B$.
    Moreover, the branch set $C \subseteq N$ is an embedded graph with no contractible components, and $\overline{V} \coloneqq p(V) \subseteq N$ is a closed handlebody that deformation retracts onto $C$.
\end{lemma}

For a natural example where $B$ is not a just disjoint union of circles, consider the action of $G = (\Z_2)^3$ on the unit 3-sphere in $\R^4$ generated by $\pi$ rotations in the 6 coordinate planes.
Then $B$ is the union of the unit circles in the each coordinate plane.

\begin{proof}[Proof of Lemma \ref{lem:branch_set_regularity}.]
    We may assume $G$ acts by isometries with respect to some Riemannian metric on $M$.
    Then for each $g \in G$, the fixed point set $F_g$ is a disjoint union of closed totally geodesic submanifolds \cite[Thm~5.1]{kobayashi}.
    By compactness, each $F_g$ has finitely many components.
    Furthermore, $F_g$ is 1-dimensional if $g$ is nontrivial; this is because for any $x \in F_g$, the element $g$ maps to a finite order element of $\Isom^+(T_xM) \cong \SO(3)$, and hence fixes a 1-dimensional subspace of $T_xM$.
    A similar local analysis shows that if $g_1$ and $g_2$ are two elements of $G$ that both fix some $x \in M$, then the components of $F_{g_1}$ and $F_{g_2}$ that contain $x$ are either equal or intersect transversely.
    Thus $B = \bigcup_{g \in G} F_g$ is a finite union of circles in $M$ intersecting transversely, and hence $B$ has the structure of an embedded $G$-invariant graph in $M$, where the vertices are the intersection points of the fixed sets $F_g$ (if $F_g$ is disjoint from all other fixed sets, we choose an arbitrary point on it to be a vertex of valence $2$).
    It follows that $B$ has no vertices of valence $1$, and hence has no contractible components.

    Now, we can build the $G$-invariant closed handlebody neighborhood $V$ of $B$.
    Choose representatives $g_1, \ldots, g_k$ of the nontrivial conjugacy classes of $G$.
    Observe that for each $g \in G$, $gF_{g_i} = F_{gg_ig^{-1}}$.
    Then each $F_{g_i}$ has a $C_{G}(g_i)$-invariant closed tubular neighborhood $N_{g_i}$ \cite[Thm~2.2]{bredon}, where $C_G(g_i)$ denotes the centralizer of $g_i$ in $G$.
    Then any $g \in G$ pushes $N_{g_i}$ to a $C_G(gg_ig^{-1})$-invariant tubular neighborhood $N_{gg_ig^{-1}}$ of $F_{gg_ig^{-1}}$ (this is well-defined, since if $hg_ih^{-1} = gg_ig^{-1}$, then $h^{-1}g \in C_G(g_i)$, and hence $gN_{g_i} = hN_{g_i}$).
    It follows that $V \coloneqq \bigcup_{i=1}^k\bigcup_{g \in G} N_{gg_ig^{-1}}$ is $G$-invariant closed neighborhood of $B$.
    By choosing the $N_{g_i}$ small enough, we can ensure that $V$ is in fact a handlebody.

    Since $B$ is a $G$-invariant embedded graph, it follows that $C$ is an embedded graph in $N$.
    We know $C$ has no contractible components, since near any $b \in B$, the stabilizer $\mathrm{Stab}_G(b) \leq G$ acts as a finite subgroup of $\SO(3)$, and thus $p(b)$ will either lie on an edge of $C$ or be a vertex of valence $\geq 2$.
    Since a quotient of a handlebody by a finite group action is still a handlebody (see e.g.\ \cite{reni-zimmermann} and/or the discussion below), it follows that $\overline{V}$ is a handlebody that deformation retracts onto $C$.
\end{proof}

\subsection{Group actions on handlebodies}\label{subsec:group_actions_on_handlebodies}
The restriction of $p$ to $V \rightarrow \overline{V}$ is precisely the quotient map for the action of $G$ on $V$.
Finite group actions on handlebodies are well-understood, as explained e.g.\ in the work of Reni-Zimmermann \cite[\S 2]{reni-zimmermann}; we will briefly summarize their discussion.

Using the Equivariant Loop Theorem \cite{meeks-yau-loop}, we can partition $V$ into $G$-orbits of $0$-handles and $1$-handles.
For each $g \in G$ and $0$-handle $E \cong B^3$, either $g(E) \cap E = \varnothing$ or $g$ acts on $E$ by an element of $\SO(3)$ (any finite group action on the $3$-ball is equivalent to a linear action \cite{kwasik-schultz-actions-on-3ball}).
For each $g \in G$ and $1$-handle $F \cong D^2 \times I$, either $g(F) \cap F = \varnothing$ or $g$ rotates the $D^2$-factor of $F$.
Note that since $V$ was constructed as a neighborhood of $B$, every $0$-handle and $1$-handle has points fixed by some $g \in G$.

For each $1$-handle $F$, the image of $F$ under $p$ is a $1$-handle $\overline{F} \cong D^2 \times I$.
The intersection $C \cap \overline{F}$ is the core arc $\{(0,0)\} \times I$.
Let $\gamma_{\overline{F}}$ denote the simple closed curve $\del D^2 \times \{\frac{1}{2}\} \subseteq \del \overline{F}$.
The curve $\gamma_{\overline{F}}$ is a meridian of $\overline{V}$, meaning it bounds a disk in $\overline{V}$; we call $\gamma_{\overline{F}}$ a \emph{branch meridian} of $p$.
The branch meridians $\gamma_{\overline{F}}$ for each such $1$-handle $\overline{F}$ normally generate the kernel of the map $\pi_1(\del \overline{V}) \rightarrow \pi_1(\overline{V})$.
If we suppose that $\gamma_{\overline{F}}$ maps to an order $k_{\overline{F}}$ element of $G$ under the monodromy homomorphism $\varphi:\pi_1(N^\circ) \rightarrow G$, then any lift of $\gamma_{\overline{F}}$ will be an arc in $\del V$ that wraps $(1/k_{\overline{F}})$ of the way around a $1$-handle in the $G$-orbit of $F$.
It follows that the lifts of the curves $\gamma_{\overline{F}}^{k_{\overline{F}}}$ for each $1$-handle $\overline{F}$ will form a system of meridians of $V$ that normally generate the kernel of the map $\pi_1(\del V) \rightarrow \pi_1(V)$.

For each $0$-handle $E$, the image of $E$ under $p$ is a $0$-handle $\overline{E}$.
Since the stabilizer $\mathrm{Stab}_G(E)$ acts on $E$ by a finite subgroup of $\SO(3)$, the quotient $\overline{E}$ inherits one of five types of orbifold structures, depending on whether $\mathrm{Stab}_G(E)$ acts via a cyclic group, dihedral group, $A_4$, $S_4$, or $A_5$.
It follows that the intersection $C \cap \overline{E}$ is either a diameter of $\overline{E}$, or a union of $3$ radii meeting at the center of $\overline{E}$ (see \cite[\S 2]{reni-zimmermann} for an illustration).

\subsection{Liftable and symmetric mapping class groups}\label{subsec:liftable_and_symmetric_mcgs}
We first define the mapping class group of the pair $(N,C)$.
Let $\Homeo^+(N,C)$ denote the group of orientation-preserving homeomorphisms of $N$ that preserve the set $C$; we do \emph{not} require that the orientation of $C$ is preserved.
Then we define the mapping class group $\Mod(N,C)$ to be $\pi_0(\Homeo^+(N,C))$, i.e.\ the quotient of $\Homeo^+(N,C)$ where two homeomorphisms are identified if they differ by an isotopy that preserves the set $C$ (we do not require isotopies to fix $C$ pointwise).
Note that any element of $\Homeo^+(N,C)$ restricts to a homeomorphism of $N - C$, and hence there is a map $\Mod(N,C) \rightarrow \Mod(N-C)$, but it is generally not an isomorphism (in particular, a homeomorphism of $N - C$ need not extend to $N$).
Moreover, since $N^\circ$ is homotopy equivalent to $N - C$, there is a natural map $\Phi:\Mod(N,C) \rightarrow \Out(\pi_1(N^\circ))$.
The map $\Phi$ is also not an isomorphism in general; in particular, any sphere twist lies in $\Ker(\Phi)$ (see e.g.\ \cite[\S 2]{hatcher-wahl}).

Now, we can define the liftable and symmetric mapping class groups.
The branched cover $p:M \rightarrow N$ restricts to an unbranched cover $M-B \rightarrow N-C$, which is classified by a surjective homomorphism $\varphi:\pi_1(N^\circ) \rightarrow G$ with $\Ker(\varphi) \cong \pi_1(M^\circ)$; we call $\varphi$ the \emph{monodromy homomorphism} of the cover $p$.
A homeomorphism $f \in \Homeo^+(N,C)$ is \emph{liftable} if the induced outer automorphism of $\pi_1(N^\circ)$ preserves the subgroup $\pi_1(M^\circ)$.
Basic algebraic topology implies that if $f$ is liftable, then the restriction $f\v_{N - C}$ lifts to a homeomorphism of $M - B$; we will show below that this lift extends to a homeomorphism of $M$ that lifts $f$.
We define the \emph{liftable mapping class group} $\LMod_p(N)$ to be the subgroup of $\Mod(N,C)$ consisting of mapping classes represented by liftable homeomorphisms.
In other words, $\LMod_p(N)$ is the stabilizer of the subgroup $\pi_1(M^\circ) \leq \pi_1(N^\circ)$. 
Since $\pi_1(M^\circ)$ has finite index in $\pi_1(N^\circ)$, the subgroup $\LMod_p(N)$ has finite index in $\Mod(N,C)$.

A homeomorphism $f$ of $M$ is \emph{symmetric} if $p(f(x)) = p(f(y))$ whenever $p(x) = p(y)$.
Any symmetric homeomorphism of $M$ descends to a homeomorphism of the pair $(N,C)$.
We define the \emph{symmetric mapping class group} $\SMod_p(M)$ to be the subgroup of the mapping class group $\Mod(M) \coloneqq \pi_0(\Homeo^+(M))$ consisting of the elements represented by symmetric homeomorphisms.
We emphasize that elements of $\SMod_p(M)$ are taken up to \emph{arbitrary} isotopy, not necessarily isotopy through symmetric homeomorphisms, so an element of $\SMod_p(M)$ does not necessarily descend to a well-defined element of $\Mod(N,C)$.
The group $\SMod_p(M)$ is contained in the normalizer of $\Gamma \leq \Mod(M)$, where $\Gamma$ is the image of $G$ in $\Mod(M)$; it would be of interest to determine whether $\SMod_p(M)$ is the full normalizer of $\Gamma$ (Birman-Hilden \cite{birman-hilden} proved that it is the full normalizer in the case of hyperbolic suraces).

\subsection{The lifting map}\label{subsec:lifting_map}
Finally, we can verify that we get a well-defined lifting homomorphism.
The main technical point to check is the following.
\begin{lemma}\label{lem:lift_exists}
    Let $f$ be a liftable homeomorphism of $N$, and let $\widetilde{f}'$ be a lift of $f\v_{N - C}$ to $M - B$.
    Then $\widetilde{f}'$ extends to a homeomorphism $\widetilde{f}$ of $M$ which lifts $f$.
\end{lemma}
\begin{proof}
    We define the extension $\widetilde{f}$ as follows.
    Fix a base point $x_0 \in M - B$, and let $y_0 = p(x_0) \in N - C$.
    For each $b \in B$, choose a path $\gamma:[0,1] \rightarrow M$ from $x_0$ to $b$, with $\gamma(t) \in M - B$ for $t \in [0,1)$.
    Then $\widetilde{f}'(\gamma\v_{[0,1)})$ is a path $[0,1) \rightarrow M - B$ based at $\widetilde{f}'(x_0)$.
    We claim that this path extends to a path $\ep:[0,1] \rightarrow M$ with $\ep(1) \in B$; assuming this claim, we define $\widetilde{f}(b) \coloneqq \ep(1)$.
    To prove the claim, let $\delta:[0,1] \rightarrow N$ be the image of $\gamma$ under $p$.
    Then $f(\delta)$ is a path in $N$ from $f(y_0)$ to $f(p(b))$.
    We can lift $f(\delta)\v_{[0,1)}$ to a path in $M - B$ based at $\widetilde{f}'(x_0)$.
    Since $p \circ \widetilde{f}' = f \circ p$ on $M -B$, this lifted path is precisely $\widetilde{f}'(\gamma\v_{[0,1)})$.
    Since $f(\delta)\v_{[0,1)}$ extends to a path ending at $f(p(b))$, it follows that $\widetilde{f}'(\gamma\v_{[0,1)})$ extends to a path ending at a point over $f(p(b))$.

    To show that this is well-defined, suppose we choose two such paths $\gamma_1$ and $\gamma_2$, and let $b_1$ and $b_2$ be the resulting two choices for $\widetilde{f}(b)$.
    Then both $b_1$ and $b_2$ lie over $f(p(b))$, and so $b_2 = g_0b_1$ for some $g_0 \in G$.
    If we assume that $G$ fixes some metric on $M$, then since $G$ is finite, the set of values $d(b_1,gb_1)$ for $g \in G - \textrm{Stab}(b_1)$ has a positive lower bound, where $d$ is the distance function of $M$.
    On the other hand, since $\gamma_1$ and $\gamma_2$ have the same endpoint and $\widetilde{f}'$ is continuous, it must be that for $t$ close to $1$, the points $\widetilde{f}'(\gamma_1\v_{[0,1)})(t)$ and $\widetilde{f}'(\gamma_2\v_{[0,1)})(t)$ are close together. 
    Since $\widetilde{f}'(\gamma_i\v_{[0,1)})(t)$ is close to $b_i$ for $t$ close to $1$, we conclude that $b_1 = b_2$.

    Thus, we indeed have a well-defined extension $\widetilde{f}$.
    The fact that $\widetilde{f}$ is continuous can be checked with sequential continuity, since for any $b \in B$, any sequence $x_i \rightarrow b$ lies on some path from $x_0$ to $b$.
    To see that $\widetilde{f}$ is a homeomorphism, let $h = f^{-1}$, and let $\widetilde{h}'$ be the lift of $h$ to $M-B$ which is inverse to $\widetilde{f}'$.
    Then the extensions $\widetilde{f}$ and $\widetilde{h}$ satisfies the equations $ \widetilde{f} \circ \widetilde{h} = \widetilde{h} \circ \widetilde{f} = id_M$ on the dense subset $M - B$, and thus the equalities hold everywhere by continuity.
    The fact that $\widetilde{f}$ lifts $f$ follows by construction.
\end{proof}

Note that since $M - B$ is dense in $M$, the extension in Lemma \ref{lem:lift_exists} is unique.
Since any lift is well-defined up to a choice of deck transformation, we get a well-defined map 
\begin{equation*}
    \widehat{\L}_p:\textrm{LHomeo}_p(N) \rightarrow \textrm{SHomeo}_p(M)/G,
\end{equation*}
where $\textrm{LHomeo}_p(N)$ is the group of liftable homeomorphisms and $\textrm{SHomeo}_p(M)$ is the group of symmetric homeomorphisms.
Moreover, Lemma \ref{lem:lift_exists} implies that for any liftable homeomorphism $f$ of $N$ with a lift $\widetilde{f}$ to $M$, any isotopy of $f$ preserving $C$ lifts to an isotopy of $\widetilde{f}$.
Thus, $\widetilde{\L}_p$ descends to a well-defined lifting map 
\begin{equation*}
    \L_p:\LMod_p(N) \rightarrow \SMod_p(M)/\Gamma,
\end{equation*}
where $\Gamma$ is the image of $G$ in $\Mod(M)$.

\begin{remark}\label{rem:BH_prop_vs_injectivity}
    As mentioned in Section \ref{subsec:intro}, we say that the cover $p$ has the \emph{Birman-Hilden property} if it satisfies the following condition: if two symmetric homeomorphisms of $M$ are isotopic, then they are isotopic through symmetric homeomorphisms.
    Observe that if $p$ has the Birman-Hilden property, then $\L_p$ is injective.
    Indeed, if $[f] \in \Ker(\L_p)$, then $f$ lifts to a symmetric homeomorphism $\widetilde{f}$ of $M$ which is isotopic to $\id_M$.
    Since $p$ has the Birman-Hilden property, there is an isotopy from $\widetilde{f}$ to $\id_M$ through symmetric homeomorphisms.
    This descends to an isotopy from $f$ to $\id_N$, which means that $[f]$ is trivial.

    The converse is slightly more subtle.
    In general, if $\L_p$ is injective, then $p$ has the following weaker property: if $f$ is a symmetric homeomorphism of $M$ which is isotopic to $\id_M$, then $f$ is isotopic to some $g \in G$ through symmetric homeomorphisms.
    This is because if $\L_p$ is injective, then any symmetric homeomorphism $f$ of $M$ which is isotopic to $\id_M$ projects to a homeomorphism $\overline{f}$ which is isotopic to $\id_N$; lifting this isotopy yields an isotopy through symmetric homeomorphisms from $f$ to a deck transformation.
    If the map $G \rightarrow \Mod(M)$ is injective (i.e.\ if $\Gamma \cong G$), then this weaker property implies the Birman-Hilden property.
    However, in general, one must be wary of the possibility that there is some $g \in G$ which is isotopic to $\id_M$ but not isotopic through symmetric homeomorphisms; determining whether this can happen appears to be a subtle problem.
\end{remark}

\section{The relative fundamental group}\label{sec:rel_fundamental_group}
Let $p:M \rightarrow N$ be a branched cover as described at the start of Section \ref{sec:preliminaries}.
In this section, we introduce the relative fundamental group $\pi_1^{\rel}(M)$, which is the key tool for proving Theorems \ref{mainthm:BH_counterexamples} and \ref{mainthm:hyperelliptic_kernel}.
In Section \ref{subsec:group_structure}, we define $\pi_1^{\rel}(M)$ and prove Proposition \ref{prop:rel_pi1_is_free_product}, which computes the group $\pi_1^{\rel}(M)$ explicitly as a quotient of $\pi_1(N^\circ)$.
In Section \ref{subsec:mcg_actions}, we show that there is a well-defined map $\Psi:\LMod_p(N) \rightarrow \Out(\pi_1^{\rel}(M))$.
In Section \ref{subsec:factoring_lifting_map}, we prove Proposition \ref{prop:lifting_map_factors}, which shows that the lifting map $\L_p$ virtually factors through $\Psi$.

\subsection{The group structure}\label{subsec:group_structure}
Let $M'$ denote $M$ or $M^\circ$.
Fix a base point $*$ on $M^\circ$, and let $G*$ denote the orbit of $*$ under $G$.
We define $\pi_1^{\rel}(M')$ to be the set $\pi_1(M', G*, *)$, i.e.\ the set of arcs from $*$ into $G*$ up to homotopy rel $G*$.
A priori, $\pi_1^{\rel}(M')$ is just a pointed set, where the distinguished element is the constant arc at $*$.
However, we can use the action of $G$ to equip $\pi_1^{\rel}(M')$ with a group structure for which the distinguished element is the identity.

First, observe that the orbit map $G \rightarrow G*$ gives us an identification $\pi_0(G*) \cong G$ mapping $*$ to $1$.
Composing this with the natural boundary map $\pi_1^{\rel}(M') \rightarrow \pi_0(G*)$, we get a function
\begin{equation*}
    \psi:\pi_1^{\rel}(M') \rightarrow G.
\end{equation*}
Now, given $[a],[b] \in \pi_1^{\rel}(M')$, we define the group operation by
\begin{equation*}
    [a][b] \coloneqq [a \cdot (\psi(a)b)],
\end{equation*}
where $\cdot$ denotes concatenation of paths.
A direct computation shows that this operation is associative.
The identity element is the constant path at $*$, and inverses are given by 
\begin{equation*}
    [a]^{-1} = [\psi(a)^{-1}\overline{a}],
\end{equation*}
where $\overline{a}$ is the reverse of $a$.

By mapping a loop in $N^\circ$ to its lift at $*$, we get a bijection $\pi_1(N^\circ) \rightarrow \pi_1^{\rel}(M^\circ)$.
With respect to the group operation above, this bijection is an isomorphism of groups.
Moreover, the natural map $\pi_1^{\rel}(M^\circ) \rightarrow \pi_1^{\rel}(M)$ is surjective, since any element of $\pi_1^{\rel}(M)$ has a representative that avoids a neighborhood of $B$.
Thus, if we equip $\pi_0(G*)$ with a group structure via the identification $\pi_0(G*) \cong G$, we get the following commutative diagram with exact rows:
\begin{equation}\label{eqn:htpy_group_diagram}
    \begin{tikzcd}
        1 & {\pi_1(M^\circ)} & {\pi_1(N^\circ)} & {G} & 1 \\
        1 & {\pi_1(M^\circ)} & {\pi_1^{\rel}(M^\circ)} & {\pi_0(G*)} & 1 \\
        1 & {\pi_1(M)} & {\pi_1^{\rel}(M)} & {\pi_0(G*)} & 1
        \arrow[from=1-1, to=1-2]
        \arrow[from=1-2, to=1-3]
        \arrow[equals, from=1-2, to=2-2]
        \arrow["\varphi", from=1-3, to=1-4]
        \arrow["\cong", from=1-3, to=2-3]
        \arrow[from=1-4, to=1-5]
        \arrow["\cong", from=1-4, to=2-4]
        \arrow[from=2-1, to=2-2]
        \arrow[from=2-2, to=2-3]
        \arrow[two heads, from=2-2, to=3-2]
        \arrow[from=2-3, to=2-4]
        \arrow[two heads, from=2-3, to=3-3]
        \arrow[from=2-4, to=2-5]
        \arrow[equals, from=2-4, to=3-4]
        \arrow[from=3-1, to=3-2]
        \arrow[from=3-2, to=3-3]
        \arrow[from=3-3, to=3-4]
        \arrow[from=3-4, to=3-5]
    \end{tikzcd}
\end{equation}
\noindent The map $\varphi$ is the monodromy homomorphism of the cover $M^\circ \rightarrow N^\circ$, and the second and third rows come from the long exact sequences of the pairs $(M^\circ, G*)$ and $(M, G*)$ respectively.

Composing the two middle vertical maps, we get a surjection $\pi_1(N^\circ) \rightarrow \pi_1^{\rel}(M)$.
The diagram implies that the kernel of this surjection is equal to the kernel of the map $\pi_1(M^\circ) \rightarrow \pi_1(M)$.
In fact, we can understand this kernel explicitly.
The inclusion $N^\circ \rightarrow N$ corresponds to filling the boundary components $\Sigma_1, \ldots, \Sigma_m$ of $N^\circ$ with handlebodies $V_1, \ldots, V_m$.
For $1 \leq i \leq m$, let $\gamma_{i1}, \ldots, \gamma_{in_i}$ be the branch merdians of $V_i$, as defined in Section \ref{subsec:group_actions_on_handlebodies}.
Let $k_{ij}$ be the order of $\varphi(\gamma_{ij}) \in G$ (the element $\varphi(\gamma_{ij})$ is well-defined up to conjugacy, so its order is well-defined).
Then we get the following.

\begin{proposition}\label{prop:rel_pi1_is_free_product}
    Let $w_{ij} \in \pi_1(N^\circ)$ be an element in the conjugacy class determined by $\gamma_{ij}$.
    Then the kernel of the map $\pi_1(N^\circ) \rightarrow \pi_1^{\rel}(M)$ is the $\pi_1(N^\circ)$-normal closure of the set of elements
    \begin{equation*}
        \{w_{ij}^{k_{ij}} \mid 1 \leq i \leq m, 1 \leq j \leq n_i \}.
    \end{equation*}
    Moreover, suppose $N^\circ \cong W_1^\circ \# \cdots \# W_r^\circ$ for some compact $3$-manifolds $W_\ell^\circ$, so there is a free splitting $\pi_1(N^\circ) \cong \pi_1(W_1^\circ) * \cdots * \pi_1(W_r^\circ)$.
    Then each $w_{ij}$ is conjugate to an element $v_{ij} \in \pi_1(W_{\ell_i}^\circ)$ for some $1 \leq \ell_i \leq r$, and the free splitting of $\pi_1(N^\circ)$ descends to a free splitting $\pi_1^{\rel}(M) \cong H_1 * \cdots * H_r$, where $H_\ell$ is the quotient of $\pi_1(W_\ell^\circ)$ by the $\pi_1(W_{\ell}^\circ)$-normal closure of the elements $v_{ij}^{k_{ij}}$ with $\ell_i = \ell$.
\end{proposition}
\begin{proof}
    As noted above, the commutative diagram (\ref{eqn:htpy_group_diagram}) implies that the kernel of the map $\pi_1(N^\circ) \rightarrow \pi_1^{\rel}(M)$ is precisely the kernel of the map $\pi_1(M^\circ) \rightarrow \pi_1(M)$.
    From Section \ref{subsec:group_actions_on_handlebodies}, this kernel is the $\pi_1(M^\circ)$-normal closure of the elements $gw_{ij}^{k_{ij}}$ for $g \in G$ (the element $gw_{ij}^{k_{ij}}$ is defined up to conjugacy, so this normal closure is well-defined).
    But by the diagram (\ref{eqn:htpy_group_diagram}), this normal closure is precisely the $\pi_1(N^\circ)$-normal closure of the elements $w_{ij}^{k_{ij}}$.
    
    For the second statement, the fact that $w_{ij}$ has such a conjugate $v_{ij}$ is because each boundary component $\Sigma_i$ lies in a single summand $W_\ell^\circ$.
    Then from the first part of the proposition, the quotient $\pi_1(N^\circ) \rightarrow \pi_1^{\rel}(M)$ is obtained starting with the free product $\pi_1(W_1^\circ) * \cdots * \pi_1(W_r^\circ)$ and then adding the relations $v_{ij}^{k_{ij}} = 1$, which yields the free product $H_1 * \cdots * H_r$.
\end{proof}

\begin{remark}\label{rem:orbifold_pi1}
    Viewing $N$ as the quotient $M/G$, we get a natural orbifold structure on $N$.
    The relative fundamental group $\pi_1^{\rel}(M)$ is then isomorphic to the orbifold fundamental group $\pi_1^{\mathrm{orb}}(N)$.
    To see this, one can observe that the presentation of $\pi_1^{\rel}(M)$ given by Proposition \ref{prop:rel_pi1_is_free_product} agrees with the usual presentation for the orbifold fundamental group (see e.g.\ \cite{haefliger-du}).
    However, we opt to use the relative fundamental group to avoid the more complicated machinery of orbifolds.
\end{remark}

\subsection{The liftable mapping class group action}\label{subsec:mcg_actions}
Next, we show that $\LMod_p(N)$ acts on $\pi_1^{\rel}(M)$ by outer automorphisms.
Let $\Phi:\LMod_p(N) \rightarrow \Out(\pi_1(N^\circ))$ be the natural map.
The definition of $\LMod_p(N)$ implies that the image of $\Phi$ lies in the subgroup $\Out(\pi_1(N^\circ), \pi_1(M^\circ))$ of outer automorphisms that preserve the normal subgroup $\pi_1(M^\circ) \leq \pi_1(N^\circ)$.
From the diagram (\ref{eqn:htpy_group_diagram}), we have a quotient map $\pi_1(N^\circ) \rightarrow \pi_1^{\rel}(M)$ taking $\pi_1(M^\circ)$ to $\pi_1(M) \leq \pi_1^{\rel}(M)$.
Then we get the following.

\begin{lemma}\label{lem:action_on_rel_pi1}
    For any $\alpha \in \LMod_p(N)$, the outer automorphism $\Phi(\alpha)$ of $\pi_1(N^\circ)$ induces an outer automorphism of $\pi_1^{\rel}(M)$ which preserves the normal subgroup $\pi_1(M) \leq \pi_1^{\rel}(M)$.
    Consequently, the map $\Phi:\LMod_p(N) \rightarrow \Out(\pi_1(N^\circ), \pi_1(M^\circ))$ descends to a well-defined map 
    \begin{equation*}
        \Psi:\LMod_p(N) \rightarrow \Out(\pi_1^{\rel}(M), \pi_1(M)).
    \end{equation*}
\end{lemma}

We note that if $f$ is a symmetric homeomorphism of $M$ which fixes the base point $* \in M$, then $f$ induces an automorphism of $\pi_1^{\rel}(M)$ which restricts to the usual induced map on $\pi_1(M)$.
From this, one can directly define the map $\Psi$ and prove that it is well-defined.
However, we opt to give a simpler proof using Proposition \ref{prop:rel_pi1_is_free_product}.
We also warn the reader that the group $\SMod_p(M)$ generally does \emph{not} act on $\pi_1^{\rel}(M)$.
This is because if $f$ is symmetric homeomorphisms, the induced map on $\pi_1^{\rel}(M)$ is only well-defined up to isotopy through symmetric homeomorphims, but elements of $\SMod_p(M)$ are taken up to arbitrary isotopy.

\begin{proof}[Proof of Lemma \ref{lem:action_on_rel_pi1}]
    Let $\gamma_{ij}$ for $1 \leq i \leq m$ and $1 \leq j \leq n_i$ be the curves in Proposition \ref{prop:rel_pi1_is_free_product}.
    By the commutative diagram (\ref{eqn:htpy_group_diagram}) and Proposition \ref{prop:rel_pi1_is_free_product}, it suffices to show that $\Phi(\alpha)$ takes the conjugacy class of each $\gamma_{ij}^{k_{ij}}$ to a conjugacy class in $\Ker(\pi_1(N^\circ) \rightarrow \pi_1^{\rel}(M))$.
    In other words, it suffices to show that if $w \in \pi_1(N^\circ)$ lies in the conjugacy class determined by $\alpha(\gamma_{ij}^{k_{ij}})$, then $w$ maps to a trivial element of $\pi_1^{\rel}(M)$.

    To see this, observe that $\alpha$ must map $\gamma_{ij}$ to a meridian $\delta$ of some boundary component $\Sigma'$ of $N^\circ$.
    Under the monodromy map $\varphi:\pi_1(N^\circ) \rightarrow G$, any element in the conjugacy class of $\gamma_{ij}$ maps to an order $k_{ij}$-element of $G$.
    Since $\Phi(\alpha)$ preserves $\pi_1(M^\circ)$, the same is true for any element in the conjugacy class of $\delta$.
    This means that each lift of $\delta^{k_{ij}}$ is a meridian of a boundary component of $M^\circ$, and hence any element in the conjugacy class of $\delta^{k_{ij}} = \alpha(\gamma_{ij}^{k_{ij}})$ maps to a trivial element of $\pi_1^{\rel}(M)$.
\end{proof}

\subsection{Factoring the lifting map}\label{subsec:factoring_lifting_map}
Finally, we can relate the map $\L_p$ to the map $\Psi$.
Let $\Gamma_*$ denote the image of $G$ in $\Out(\pi_1(M))$, and let $\Out_{\Gamma_*}(\pi_1(M))$ denote the normalizer of $\Gamma_*$.
Then we have a natural map
\begin{equation*}
    \Theta: \SMod_p(M)/\Gamma \rightarrow \Out_{\Gamma_*}(\pi_1(M))/\Gamma.
\end{equation*}
Our goal is to show that the composition $\Theta \circ \L_p$ factors through $\Psi$.

To this end, we observe first that there is a well-defined map 
\begin{equation*}
    \overline{r}:\Out(\pi_1^{\rel}(M), \pi_1(M)) \rightarrow \Out_{\Gamma_*}(\pi_1(M))/\Gamma_*.
\end{equation*}
To obtain $\overline{r}$, we start with the natural restriction map 
\begin{equation*}
    r:\Aut(\pi_1^{\rel}(M), \pi_1(M)) \rightarrow \Aut(\pi_1(M)).
\end{equation*}
Post-composing with the projection $\Aut(\pi_1(M)) \rightarrow \Out(\pi_1(M))$, we get a map 
\begin{equation*}
    r':\Aut(\pi_1^{\rel}(M), \pi_1(M)) \rightarrow \Out(\pi_1(M)).
\end{equation*}
From the short exact sequence $\pi_1(M) \hookrightarrow \pi_1^{\rel}(M) \twoheadrightarrow G$, we see that $r'$ takes $\Inn(\pi_1^{\rel}(M))$ to the subgroup $\Gamma_* \leq \Out(\pi_1(M))$.
Since $\Inn(\pi_1^{\rel}(M))$ is normal in $\Aut(\pi_1^{\rel}(M), \pi_1(M))$, it follows that $\Gamma_*$ is normal in the image of $r'$, i.e.\ the image of $r'$ lies in $\Out_{\Gamma_*}(\pi_1(M))$.
Then $r'$ descends to the map $\overline{r}$.

\begin{proposition}\label{prop:lifting_map_factors}
    We have the following commutative diagram:
    \begin{equation*}
    \begin{tikzcd}
        {\LMod_p(N)} & {\Out(\pi_1^{\rel}(M), \pi_1(M))} \\
        {\SMod_p(M)/\Gamma} & {\Out_{\Gamma_*}(\pi_1(M))/\Gamma_*}
        \arrow["\Psi", from=1-1, to=1-2]
        \arrow["{\L_p}"', from=1-1, to=2-1]
        \arrow["{\overline{r}}", from=1-2, to=2-2]
        \arrow["\Theta"', from=2-1, to=2-2]
    \end{tikzcd}
    \end{equation*}
    Moreover, the kernel of the map $\Theta$ is finite, and consequently if $\alpha \in \Ker(\Psi)$, then $\L_p(\alpha)$ has finite order.
\end{proposition}
\begin{proof}
    First, we can prove the commutativity of the diagram.
    Recall that the composition $\Theta \circ \L_p$ acts as follows: given $\alpha \in \LMod_p(N)$, choose a representative $f$, then choose a lift $\widetilde{f}$ to $M$, and then take the outer automorphism $\omega$ of $\pi_1(M)$ induced by $\widetilde{f}$.
    On the other hand, the map $\widetilde{f}$ also induces an outer automorphism $\widehat{\omega}$ of $\pi_1(M^\circ)$, and $\widehat{\omega}$ is precisely the restriction of $\Phi(\alpha) \in \Out(\pi_1(N^\circ), \pi_1(M^\circ))$.
    By the commutative diagram (\ref{eqn:htpy_group_diagram}) and Lemma \ref{lem:action_on_rel_pi1}, $\Phi(\alpha)$ descends to the outer automorphism $\Psi(\alpha)$ of $\pi_1^{\rel}(M)$, and its restriction $\widehat{\omega}$ descends to $\omega$.

    For the second statement, suppose $\widetilde{f}$ is a symmetric homeomorphism of $M$ which induces the trivial outer automorphism of $\pi_1(M)$.
    Then $\widetilde{f}$ is isotopic to a homeomorphism which fixes the base point $* \in M^\circ$, and acts trivially on $\pi_1(M)$.
    Since $M$ is closed, it follows from \cite[Prop~2.1]{hatcher-wahl} that $\widetilde{f}$ is isotopic to a product of sphere twists.
    The subgroup of $\Mod(M)$ generated by sphere twists is finite \cite[Prop~1.2]{mccullough}, and thus $\Ker(\Theta)$ is finite.
\end{proof}

\section{Constructing elements in the lifting kernel}\label{sec:pf_thm_A}

In this section, we prove Theorem \ref{mainthm:BH_counterexamples} and Corollary \ref{maincor:BH_counterexamples}.
The idea is to find certain homeomorphisms of $N^\circ$, called \emph{slide homeomorphisms}, which induce an infinite order automorphism of $\pi_1(N^\circ)$ but act trivially on $\pi_1^{\rel}(M)$.

Before we begin, we fix some notation for this section.
Let $p:M \rightarrow N$ be a branched cover as described at the start of Section \ref{sec:preliminaries}.
Throughout this section, we suppose that $N^\circ$ is homeomorphic to a connect sum $W_1^\circ \# W_2^\circ$ for some compact orientable $3$-manifolds $W_i^\circ$ with no spherical boundary components.
We fix an embedded separating sphere $S \subseteq N^\circ$ which realizes this decomposition, and we let $N_i^\circ \subseteq N^\circ$ be the closure of the component of $N^\circ - S$ corresponding to $W_i^\circ$ (so $W_i^\circ$ is obtained by capping the spherical boundary component of $N_i^\circ$).
We fix a base point $*$ on $S \subseteq N^\circ$, and a basepoint $* \in M^\circ$ above it.
We then get a natural free splitting 
\begin{equation*}
    \pi_1(N^\circ) \cong \pi_1(N_1^\circ) * \pi_1(N_2^\circ) \cong \pi_1(W_1^\circ) * \pi_1(W_2^\circ).
\end{equation*}
We will mostly be interested in the case that $N^\circ$ is not prime, and the manifolds $W_i^\circ$ are not $3$-spheres (and hence $\pi_1(W_i^\circ)$ is nontrivial by the Poincar\'e conjecture).

In Section \ref{subsec:slide_homeos}, we define slide homeomorphisms, and show that certain slide homeomorphisms have infinite order.
In Section \ref{subsec:pf_thm_A_case_(i)}, we prove part (i) of Theorem \ref{mainthm:BH_counterexamples} by finding infinite order slide homeomorphisms that lift to finite order elements.
In Section \ref{subsec:pf_thm_A_case_ii}, we prove part (ii) of Theorem \ref{mainthm:BH_counterexamples}; the argument for part (i) does not apply in this case for certain Seifert fibered summands, so we tackle these cases directly.
In Section \ref{subsec:pf_cor_B}, we deduce Corollary \ref{maincor:BH_counterexamples} using the Equivariant Sphere Theorem.
Finally, in Section \ref{subsec:exceptional_cases}, we study the exceptional cases of Theorem \ref{mainthm:BH_counterexamples}.

\subsection{Slide homeomorphisms}\label{subsec:slide_homeos}
We define the notion of a slide homeomorphism following McCullough \cite[\S 1]{mccullough}, though we will take a slightly more formal approach.
Recall that the base point $* \in N^\circ$ lies on $S$.
Let $\gamma$ be a simple closed curve in $N_i^\circ$ based at $*$.
View $W^\circ_i$ as obtained from capping the boundary sphere $S$ of $N_i^\circ$ by a 3-ball $E$.

First, we introduce the notion of a boundary slide of $N_i^\circ$ along $\gamma$.
Informally, let $h_t$ be an ambient isotopy of $W_i^\circ$ with $h_0 = \id$ and $h_1\v_E = \id_E$ which moves $E$ along $\gamma$; a boundary slide is the restriction of $h_1$ to $N_i^\circ$ for any choice of such an isotopy $h_t$.
We formalize this idea as follows (for convenience, we work in the smooth category).
Choose an arc in $E$ from the center of $E$ to $*$, and conjugate $\gamma$ by this arc to obtain a loop $\gamma'$ based at the center of $E$.
Then choose a lift of $\gamma'$ to a loop $\widetilde{\gamma}'$ in the (oriented) frame bundle $\Fr^+(W_i^\circ)$.
The frame bundle $\Fr^+(W_i^\circ)$ is weakly homotopy equivalent to the space $\Emb^+(E, W_i^\circ)$ of smooth orientation-preserving embeddings of $E$ in $W_i^\circ$ \cite[Thm~2.6.C]{ivanov}.
Moreover, there is a fibration 
\begin{equation*}
    \Diff_{\del}^+(W_i^\circ \rel E) \rightarrow \Diff_{\del}^+(W_i^\circ) \rightarrow \Emb^+(E, W_i^\circ),
\end{equation*}
where $\Diff_{\del}^+(W_i^\circ)$ is the group of orientation-preserving diffeomorphisms fixing $\del W_i^\circ$ pointwise, and $\Diff_{\del}^+(W_i^\circ \rel E)$ is the group of orientation-preserving diffeomorphisms that fix $\del W_i^\circ$ and $E$ pointwise \cite[Thm~2.6.A]{ivanov}.
We define a \emph{boundary slide} of $N_i^\circ$ along $\gamma$ to be a homeomorphism which represents the image of $\widetilde{\gamma}'$ under the composition 
\begin{equation*}
    \pi_1(\Fr^+(W_i^\circ)) \xrightarrow{\cong} \pi_1(\Emb^+(E, W_i^\circ)) \rightarrow \pi_0(\Diff_{\del}^+(W_i^\circ \rel E)) \rightarrow \pi_0(\Diff_{\del}^+(N_i^\circ)),
\end{equation*}
where the second map comes from the long exact sequence on homotopy groups and the third map is induced by restriction.

The isotopy class (rel boundary) of a boundary slide depends only on the based homotopy class of the chosen lift $\widetilde{\gamma}'$.
Note that since any compact oriented $3$-manifold is parallelizable, the frame bundle is a trivial $\GL_3^+(\R)$-bundle, and hence
\begin{equation*}
    \pi_1(\Fr^+(W_i^\circ)) \cong \pi_1(\GL_3^+(\R)) \times \pi_1(W_i^\circ) \cong \Z/2\Z \times \pi_1(W_i^\circ).
\end{equation*}
Boundary slides coming the $\pi_1(W_i^\circ)$-factor behave like ``point pushing'' maps on surfaces (see \cite[\S 4.2]{farb-margalit}), while a boundary slide coming from the $\pi_1(\GL_3^+(\R))$-factor is a sphere twist about $S$ (cf.\ \cite[Rem~2.4]{hatcher-wahl}).
The based homotopy class of our lift $\widetilde{\gamma}'$ is determined by the based homotopy class of $\gamma$ up to the $\pi_1(\GL_3^+(\R))$-factor.
It follows that given an element of $\pi_1(N_i^\circ, *)$, the isotopy class of the resulting boundary slide is well-defined up to a sphere twist about $S$.

Now, we define a \emph{slide along $\gamma$} to be a homeomorphism of $N^\circ$ which acts on $N_i^\circ$ by a boundary slide along $\gamma$, and which fixes the other summand $N_j^\circ$ pointwise.
We again emphasize that given an element of $\pi_1(N_i^\circ, *)$, the isotopy class of a resulting slide homeomorphism is only well-defined up to a sphere twist about $S$.

Next, we can examine the action of a slide homeomorphism on $\pi_1(N^\circ)$.
Suppose that $\gamma$ is a simple closed curve in $N_1^\circ$ based at $*$ and $f$ is a slide along $\gamma$.
Since $f$ fixes $*$, it induces an automorphism of $\pi_1(N^\circ)$.
In particular, we see that $f$ acts on $\pi_1(N^\circ)$ by \emph{partial conjugation}.
That is, with respect to the free splitting $\pi_1(N^\circ) = \pi_1(N_1^\circ) * \pi_1(N_2^\circ)$, the homeomorphism $f$ conjugates $\pi_1(N_1^\circ)$ by $[\gamma]$ and fixes $\pi_1(N_2^\circ)$ pointwise.
The analogous statement holds if $\gamma$ is a simple closed curve in $N_2^\circ$.
This observation has the following immediate consequence.

\begin{lemma}\label{lem:finite_order_iff_central_power}
    Assume $\pi_1(N_1^\circ)$ and $\pi_1(N_2^\circ)$ are both nontrivial.
    Suppose $\gamma$ is a simple closed curve in $N_i^\circ$ based at $*$.
    Let $\alpha = [f] \in \Mod(N,C)$ where $f$ is a slide along $\gamma$.
    Then $\alpha$ has infinite order if no power of $[\gamma] \in \pi_1(N_i^\circ)$ is central.
\end{lemma}

\begin{remark}\label{rem:surface-slides}
    Slide homeomorphisms have an analogous construction on surfaces.
    Let $\Sigma$ be a compact orientable surface, and let $\delta$ be a separating simple closed curve on $\Sigma$.
    Let $\Sigma_1$ and $\Sigma_2$ denote the closures of the components of $\Sigma - \delta$, and assume each $\pi_1(\Sigma_i)$ is nontrivial.
    Let $\widehat{\Sigma}_i$ denote the surface obtained from $\Sigma_i$ by capping the boundary component coming from $\delta$ with a disk.
    We define a \emph{boundary slide} of $\Sigma_i$ to be a homeomorphism of $\Sigma_i$ which represents an element of the kernel of the map $\Mod(\Sigma_i) \rightarrow \Mod(\widehat{\Sigma}_i)$.
    Boundary slides on $\Sigma_i$ correspond to elements of the fundamental group of the unit tangent bundle $\UT(\widehat{\Sigma}_i)$ \cite[\S 4.2.5]{farb-margalit}, just as boundary slides on $N_i^\circ$ correspond to elements of $\pi_1(\Fr^+(W_i^\circ))$.
    Then, we can define a \emph{slide homeomorphism} of $\Sigma$ to be a homeomorphism which acts on $\Sigma_i$ by a boundary slide and which fixes the other subsurface $\Sigma_j$ pointwise.

    For a surface $\Sigma$, the action of a slide homeomorphism on $\pi_1(\Sigma)$ is more complicated than in the $3$-manifold case.
    This is because $\pi_1(\UT(\widehat{\Sigma}_i))$ fits into a short exact sequence 
    \begin{equation*}
        1 \rightarrow \Z \rightarrow \pi_1(\UT(\widehat{\Sigma}_i)) \rightarrow \pi_1(\widehat{\Sigma}_i) \rightarrow 1.
    \end{equation*}
    Given a generator of the kernel $\Z$, the corresponding slide homeomorphism is a Dehn twist around $\delta$, which will act nontrivially on $\pi_1(\Sigma)$, unlike sphere twists in the $3$-manifold case.
    This short exact sequence also implies that an element of $\pi_1(\Sigma_i)$ corresponds to infinitely many isotopy classes of slide homeomorphisms, unlike in the $3$-manifold case where there are at most $2$ isotopy classes.
    This second difference is ultimately a consequence of the fact that $\pi_1(\GL_2(\R)) \cong \Z$ while $\pi_1(\GL_3(\R)) \cong \Z/2\Z$.
    
\end{remark}

\subsection{Part (i) of Theorem \ref{mainthm:BH_counterexamples}}\label{subsec:pf_thm_A_case_(i)}
Now, we can prove part (i) Theorem \ref{mainthm:BH_counterexamples}.
Our goal is to construct an infinite order mapping class $\alpha \in \LMod_p(N)$ such that $\L_p(\alpha)$ has finite order.
The mapping class $\alpha$ will be represented by a slide homeomorphism around the boundary of $N^\circ$.

Recall from Section \ref{subsec:group_actions_on_handlebodies} that a \emph{branch meridian} is a simple closed curve in $\del N^\circ$ that bounds a disk in $N$ that intersects $C$ once.
We claim that if $w \in \pi_1(N^\circ)$ is freely homotopic to a branch meridian, then $w$ has infinite order.
To see this, we refine the decomposition $N^\circ = W_1^\circ \# W_2^\circ$ to a prime decomposition $N^\circ = P_1^\circ \# \cdots \# P_k^\circ$.
Since any boundary component lies on some $P_i^\circ$, the element $w$ is conjugate into the free factor $\pi_1(P_i^\circ)$ of $\pi_1(N^\circ)$.
The group $\pi_1(P_i^\circ)$ is torsion free; this is because $P_i^\circ$ is irreducible (it has nonempty boundary so it is not $S^1 \times S^2$), and any boundary component has positive genus, so the ``Half Lives, Half Dies" Theorem implies that $\pi_1(P_i^\circ)$ is infinite.
Since $w$ is nontrivial (it maps to a nontrivial deck transformation), it therefore must have infinite order.

On the other hand, it follows from the discussion in Section \ref{subsec:group_actions_on_handlebodies} that the image of $w$ in $\pi_1^{\rel}(M)$ has finite order.
From these observations, we can show that a slide along any branch meridian lies in $\Ker(\L_p)$.

\begin{lemma}\label{lem:slide_lift_finite_order}
    Let $\gamma$ be a simple closed curve in $N^\circ$ based at $*$ which is freely homotopic to a branch meridian, and let $\alpha = [f] \in \LMod_p(N)$ where $f$ is a slide along $\gamma$.
    Then $\L_p(\alpha)$ has finite order.
\end{lemma}
\begin{proof}
    Without loss of generality, assume that $\gamma$ lies in $N_1^\circ$.
    By definition, $f$ fixes the basepoint $* \in N^\circ$; let $\widetilde{f}$ be the lift fixing the base point $* \in M^\circ$.
    Recall from Proposition \ref{prop:rel_pi1_is_free_product} that the quotient map $\pi_1(N^\circ) \rightarrow \pi_1^{\rel}(M)$ preserves the free splitting $\pi_1(N^\circ) \cong \pi_1(N_1^\circ) * \pi_1(N_2^\circ)$, i.e.\ $\pi_1^{\rel}(M) = H_1 * H_2$ where $H_i$ is a quotient of $\pi_1(N_i^\circ)$.
    Let $w \in \pi_1^{\rel}(M)$ be the image of $[\gamma] \in \pi_1(N^\circ)$.
    Recall that $f$ acts on $\pi_1(N^\circ)$ by conjugating $\pi_1(N_1^\circ)$ by $[\gamma]$ and fixing $\pi_1(N_2^\circ)$ pointwise.
    It follows that $\widetilde{f}$ conjugates $H_1$ by $w$ and fixes $H_2$ pointwise.
    Since $w$ has finite order, $\widetilde{f}$ induces a finite order automorphism of $\pi_1^{\rel}(M)$, and hence by Proposition \ref{prop:lifting_map_factors}, we conclude that $\L_p(\alpha)$ has finite order.
\end{proof}

Now, we prove part (i) of \ref{mainthm:BH_counterexamples} by using Lemmas \ref{lem:finite_order_iff_central_power} and \ref{lem:slide_lift_finite_order} to find an infinite order element of $\Ker(\L_p)$.

\begin{proposition}[Theorem \ref{mainthm:BH_counterexamples} part (i)]\label{prop:BH_counterexamples_part_i}
    If $N^\circ$ has at least $3$ prime factors, then $\Ker(\L_p)$ is infinite.
\end{proposition}
\begin{proof}
    Let $N^\circ = P_1^\circ \# \cdots \# P_k^\circ$ be the prime decomposition of $N^\circ$.
    Assume without loss of generality that $P_1^\circ$ has nonempty boundary.
    Set $W_1^\circ = P_1^\circ \# P_2^\circ$, and let $W_2^\circ$ be the connected sum of the remaining prime factors, so $N^\circ \cong W_1^\circ \# W_2^\circ$.
    Since $k \geq 3$, each $\pi_1(W_i^\circ)$ is nontrivial.
    Choose a simple closed curve $\gamma$ in $N_1^\circ$ based at $*$ which is freely homotopic to a branch meridian of $N_1^\circ$.
    Recall from above that $[\gamma]$ is an infinite order element of $\pi_1(N^\circ)$.
    Since $\pi_1(N_1^\circ)$ is the nontrivial free product $\pi_1(P_1^\circ)*\pi_1(P_2^\circ)$, we know that $[\gamma] \in \pi_1(N_1^\circ) \leq \pi_1(N^\circ)$ has no power which is central in $\pi_1(N_1^\circ)$.
    So, let $\alpha = [f] \in \Mod(N,C)$ where $f$ is a slide along $\gamma$.
    Since $\LMod_p(N)$ has finite index, there is some power $\alpha^m$ which lies in $\LMod_p(N)$.
    By Lemma \ref{lem:finite_order_iff_central_power}, $\alpha$ has infinite order, and hence $\alpha^m$ does too.
    On the other hand, Lemma \ref{lem:slide_lift_finite_order} implies that $\L_p(\alpha^m)$ has finite order.
    Thus, some higher power of $\alpha$ will be an infinite order element in $\Ker(\L_p)$.
\end{proof}

\subsection{Part (ii) of Theorem \ref{mainthm:BH_counterexamples}}\label{subsec:pf_thm_A_case_ii}

It remains to study the case that $N^\circ$ has two prime factors, i.e.\ $W_1^\circ$ and $W_2^\circ$ are both prime.
If $N_i^\circ$ has a branch meridian with no power which is central in $\pi_1(N_i^\circ)$, then we can repeat the proof of Proposition \ref{prop:BH_counterexamples_part_i} to show that $\Ker(\L_p)$ is infinite.
However, it's possible that $N_i^\circ$ has an infinite center, and that a branch meridian has a power in this center.
In this case, the Seifert Fiber Space Theorem implies that $W_i^\circ$ is Seifert fibered (since $W_i^\circ$ is orientable and Haken, this follows from the work of Waldhausen \cite{waldhausen-SFST}).
Most of our attention will therefore be spent on the case that some $W_i^\circ$ is Seifert fibered.

Our first step is to show that $\Ker(\L_p)$ is infinite for certain branched covers where $N^\circ$ is Seifert fibered.

\begin{lemma}\label{lem:SF_kernel_infinite}
    Suppose that $\pi_1(N^\circ)$ has an infinite center, and that there is a Seifert fibration $N^\circ \rightarrow \Sigma$, where $\Sigma$ is a compact orientable $2$-orbifold of genus $g$ with $b$ boundary components.
    Assume that $g \geq 1$ or $b \geq 2$.
    If each branch meridian of $N^\circ$ is homotopic to a regular fiber, then there is an infinite order mapping class $\alpha \in \LMod_p(N)$ such that $\L_p(\alpha)$ has finite order.
    In particular, $\Ker(\L_p)$ is infinite.
\end{lemma}
\begin{proof}
    Let $(p_1,q_1), \ldots, (p_r,q_r)$ be the signatures of the singular fibers of $N^\circ$.
    Recall that $\pi_1(N^\circ)$ has a presentation with generators $a_1,b_1, \ldots, a_g,b_g, c_1, \ldots, c_b, d_1, \ldots, d_r$, and $z$, subject to the relations
    \begin{align*}
        &\bullet z \text{ is central}, \\
        &\bullet \prod_{i=1}^g [a_i,b_i] \prod_{i=1}^b c_i \prod_{i=1}^r d_i = z^\ell, \\
        &\bullet d_i^{p_i}z^{q_i} = 1,
    \end{align*}
    for some $\ell \in \Z$.
    Here $z$ corresponds to the regular fiber.
    We begin by defining a certain automorphism $\psi$ of $\pi_1(N^\circ)$.
    Namely, if $g \geq 1$, we define $\psi$ by mapping $a_1$ to $a_1z$, and fixing all other generators.
    Otherwise, if $b \geq 2$, we define $\psi$ by mapping $c_1$ to $c_1z$, mapping $c_2$ to $c_2z^{-1}$, and fixing all other generators.
    The relations above tell us that in both cases, the map $\psi$ is well-defined.

    Now, we claim that there is a homeomorphism $f$ of the pair $(N,C)$ which induces the outer automorphism of $\pi_1(N^\circ)$ represented by $\psi$.
    To see this, we observe that the manifold $N^\circ$ is Haken (it is orientable and has a nonspherical boundary component), each boundary component of $N^\circ$ is incompressible (this follows from our assumption on $g$ and $b$), and the automorphism $\psi$ preserves the peripheral subgroups $\la c_i, z \ra$.
    It therefore follows from a theorem of Waldhausen \cite*[Cor~6.5]{waldhauden-pi1-iso} that there is a homeomorphism $f^\circ$ of $N^\circ$ which induces $\psi$.
    So, it just remains to check that $f^\circ$ extends to a homeomorphism $f$ of the pair $(N,C)$.
    Note that since $N^\circ$ is Seifert fibered, each component of $\del N^\circ$ is a torus.
    Observe that for each boundary component $T \subseteq \del N^\circ$, the inclusion $T \hookrightarrow N^\circ$ is $\pi_1$-injective, and we can identify $\pi_1(T)$ with the subgroup $\la c_i, z \ra$ for some $i$.
    Let $U_T \subseteq N$ denote the solid torus which fills $T$ under the inclusion $N^\circ \hookrightarrow N$, so the meridian of $U_T$ is a regular fiber.
    Then from the definition of $\psi$, the restriction $f^\circ\v_T$ is either isotopic to $\id_T$, or is a Dehn twist which preserves the meridian of $U_T$.
    Thus we see that $f^\circ$ extends to a homeomorphism $f$ of $N$, and since $f$ preserves the core curve of each solid torus $U_T$ up to isotopy, we me assume up to isotopy that $f$ preserves $C$.

    Now, let $\alpha_0 = [f] \in \Mod(N, C)$, and let $\alpha$ be a power of $\alpha_0$ which lies in $\LMod_p(N)$.
    We claim that $\alpha$ has infinite order and that $\L_p(\alpha)$ has finite order; the lemma follows from this claim.
    The element $\alpha$ has infinite order since the outer automorphism $[\psi]$ does.
    To show that $\L_p(\alpha)$ has finite order, Proposition \ref{prop:lifting_map_factors} tells us that it's enough to check that $\psi$ induces a finite order automorphism of $\pi_1^{\rel}(M)$.
    But this follows because each branch meridian of $N^\circ$ has a power in the kernel of $\pi_1(N^\circ) \rightarrow \pi_1^{\rel}(M)$, which means that $z$ has a power $z^k$ in the kernel of $\pi_1(N^\circ) \rightarrow \pi_1^{\rel}(M)$, and hence $\psi^k$ induces a trivial automorphism of $\pi_1^{\rel}(M)$
\end{proof}

Lemma \ref{lem:SF_kernel_infinite} tells us that $\Ker(\L_p)$ is infinite for certain covers where $N^\circ$ is Seifert fibered.
The next step is to check that this implies $\Ker(\L_p)$ is infinite for certain covers where $N^\circ$ is reducible.

First, we expand our notation.
Note that the inclusion $N^\circ \hookrightarrow N$ induces a way to fill in each non-spherical boundary component of $N_i^\circ$; we call the resulting manifold $N_i$.
Similarly, the inclusion $N^\circ \hookrightarrow N$ induces a way to fill in each boundary component of $W_i^\circ$; we call the resulting manifold $W_i$.
So, $N_i$ has a single boundary component, namely the sphere $S$, and capping this boundary yields the closed manifold $W_i$.

\begin{lemma}\label{lem:summand_homeo_in_kernel}
    Assume $\pi_1(N_1^\circ)$ and $\pi_1(N_2^\circ)$ are both nontrivial.
    Choose a component $\widetilde{N}_i$ of $p^{-1}(N_i)$, so $p$ restricts a branched cover $q:\widetilde{N}_i \rightarrow N_i$.
    Let $q':\widetilde{W}_i \rightarrow W_i$ be the branched cover of closed manifolds obtained by capping all spherical boundary components.
    Suppose there is an infinite order element $\alpha' \in \LMod_{q'}(W_i)$ such that $\L_{q'}(\alpha')$ is trivial.
    Then there is an infinite order element $\alpha \in \LMod_p(N)$ such that $\L_p(\alpha)$ is trivial.
\end{lemma}
\begin{proof}
    View $W_i$ as being obtained from $N_i$ by capping the spherical boundary component $S$ with a 3-ball $E$.
    Choose a representative $f'$ of $\alpha'$ which fixes $E$, and hence fixes the base point $* \in N^\circ$ in particular.
    Then $f'$ induces an automorphism of $\pi_1(N_i^\circ)$ which descends to an infinite order outer automorphism and which perserves the subgroup $\pi_1(\widetilde{N}_i^\circ)$ and induces an inner automorphism on the quotient $\pi_1(\widetilde{N}_i)$.
    Post-composing with a slide homeomorphism, we can obtain a homeomorphism $f''$ which acts on $\pi_1(N_i^\circ)$ with infinite order and induces the identity map on $\pi_1(\widetilde{N}_i)$.
    Now, extend $f''$ to a homeomorphism $f$ of $(N,C)$ which acts by the identity on $N - N_i^\circ$.
    Then $f$ induces an infinite order automorphism of the free factor $\pi_1(N_i^\circ) \leq \pi_1(N^\circ)$, and hence $f$ represents an infinite order element $\alpha = [f] \in \LMod_p(N)$.
    Moreover, if we let $\widetilde{f} \in \Homeo^+(M)$ be the lift which fixes the base point $* \in M$, then by Proposition \ref{prop:rel_pi1_is_free_product} $\widetilde{f}$ induces the identity automorphism of $\pi_1^{\rel}(M)$.
    By Proposition \ref{prop:lifting_map_factors}, we conclude that $\L_p(\alpha)$ has finite order.
\end{proof}

Now, combining Lemmas \ref{lem:SF_kernel_infinite} and \ref{lem:summand_homeo_in_kernel}, we can prove part (ii) of Theorem \ref{mainthm:BH_counterexamples}.

\begin{proposition}[Theorem \ref{mainthm:BH_counterexamples} part (ii)]
    Suppose that $W_1^\circ$ and $W_2^\circ$ are both nontrivial and prime, and that $W_i^\circ$ has nonempty boundary for some $i \in \{1,2\}$.
    Suppose $W_i^\circ$ is not a Seifert fibration over an orientable compact $2$-orbifold with genus $0$ and $1$ boundary component.
    Then $\Ker(\L_p)$ is infinite.
\end{proposition}
\begin{proof}
    First, suppose that $N_i^\circ$ has a curve $\gamma$ based at $*$ which is freely homotopic to a branch meridian, and that the element $[\gamma]$ has no central power in $\pi_1(N_i^\circ)$.
    Then, as in the proof of Proposition \ref{prop:BH_counterexamples_part_i}, we deduce from Lemmas \ref{lem:finite_order_iff_central_power} and \ref{lem:slide_lift_finite_order} that $\Ker(\L_p)$ is infinite.

    Otherwise, suppose that for any based curve $\gamma$ in $N_i^\circ$ which is freely homotopic to a branch meridian, the element $[\gamma]$ has a power which is central in $\pi_1(N_i^\circ)$; let $[\gamma]^k$ be the smallest such power.
    Recall from above that $[\gamma]$ has infinite order.
    Then the Seifert Fiber Space Theorem \cite{waldhausen-SFST} says that there is a Seifert fibering of $W_i^\circ$, and $[\gamma]^k = [\delta]^\ell$ for some $\ell \geq 1$, where $\delta$ is freely homotopic to a regular fiber.
    Note that since $\pi_1(W_i^\circ)$ has an infinite center, the base $\Sigma$ of the Seifert fibering will be orientable.
    Each boundary component of $W_i^\circ$ will be a torus; by assumption $W_i^\circ$ is not a solid torus, so $W_i^\circ$ has incompressible boundary.

    We claim that in fact $[\gamma] = [\delta]$, i.e.\ $\gamma$ is freely homotopic to a regular fiber.
    Let $T \subseteq \del W_i^\circ$ be the component for which $\gamma$ is a branch meridian.
    Let $Q \leq \pi_1(W^\circ)$ be the image of $\pi_1(T) \hookrightarrow \pi_1(W_i^\circ)$.
    Choose an element $[\ep] \in Q$ such that $[\delta]$ and $[\ep]$ form a basis of $W \cong \pi_1(T)$.
    Then since $\gamma$ is freely homotopic into $T$, $[\gamma]$ is conjugate to an element of the form $[\delta]^p[\ep]^q$.
    Furthermore, $[\gamma]^k$ is conjugate to the element $[\delta]^{pk}[\ep]^{qk}$, and since $[\gamma]^k$ is central, we in fact have that $[\gamma]^k = [\delta]^{pk}[\ep]^{qk}$.
    But then 
    \begin{equation*}
        [\delta]^\ell = [\gamma]^k = [\delta]^{pk}[\ep]^{qk}.
    \end{equation*}
    This is possible only if $q = 0$, i.e.\ $[\gamma]$ is conjugate to $[\delta]^p$, and since $[\delta]$ is central, this means $[\gamma] = [\delta]^p$.
    Finally, since $\gamma$ is freely homotopic to a simple closed curve on $T$, it must be that $[\gamma]$ is conjugate to a primitive element of $W$, and so in fact $p=1$, i.e.\ $[\gamma] = [\delta]$.

    Choose a component $\widetilde{N}_i$ of $p^{-1}(N_i)$, so $p$ restricts a branched cover $q:\widetilde{N}_i \rightarrow N_i$.
    If we cap off the spherical boundary components, we get a branched cover $q':\widetilde{W}_i \rightarrow W_i$ which satisfies the hypotheses of Lemma \ref{lem:SF_kernel_infinite}.
    Thus there is an infinite order mapping class $\alpha' \in \LMod_{q'}(W_i)$ such that $\L_{q'}(\alpha)$ has finite order.
    Then by Lemma \ref{lem:summand_homeo_in_kernel}, we conclude that $\Ker(\L_p)$ is infinite.
\end{proof}

\subsection{Proof of Corollary \ref{maincor:BH_counterexamples}}\label{subsec:pf_cor_B}
Next, we can deduce Corollary \ref{maincor:BH_counterexamples} from part (i) of Theorem \ref{mainthm:BH_counterexamples} using the Equivariant Sphere Theorem of Meeks-Yau \cite{meeks-yau-sphere}.
It suffices to show the following.

\begin{lemma}
    Suppose that $G$ does not act freely on $M$, and that $\pi_2(M^\circ)_G$ has rank at least $3$ as a $\pi_1(M^\circ)$-module.
    Then $N^\circ$ has at least $3$ prime factors.
\end{lemma}
\begin{proof}
    Since $\pi_2(M^\circ)$ has rank at least $3$ as a $\pi_1(M^\circ)$-module, it follows from of Assertion 1 in the proof of \cite[Thm~9]{meeks-yau-sphere} that $M$ has $G$-invariant embedded multispheres $\calS_1$, $\calS_2$, $\calS_3$ whose components are pairwise disjoint and non-homotopic.
    This implies in particular that $N^\circ$ is not prime.
    Indeed, each $\calS_i$ descends to either an essential embedded sphere or a $1$-sided projective plane in $N^\circ$.
    If $N^\circ$ has an essential embedded sphere, then it is reducible, and hence prime since $\del N^\circ \neq \varnothing$.
    If $N^\circ$ has a $1$-sided projective plane, then it contains an $\RP^3$-summand, and again must therefore be non-prime since $\del N^\circ \neq \varnothing$.
    Thus, $N^\circ$ contains an embedded separating sphere $S_1$; without loss of generality, we may assume $\calS_1$ projects to $S_1$.

    Let $N_1^\circ$ and $N_2^\circ$ be the components of $N^\circ - S_1$.
    To prove the lemma, it's enough to show that $N_1^\circ$ and $N_2^\circ$ cannot both be prime.
    We can apply a similar argument.
    Assume without loss of generality that $N_1^\circ$ has a non-spherical boundary component.
    Let $W_i^\circ$ be obtained by capping the spherical boundary component of $N_i^\circ$.
    For $i = 2$ or $3$, the image $S_i$ of $\calS_i$ in $N^\circ$ is a nontrivial embedded sphere or 1-sided projective plane.
    If either $S_2$ or $S_3$ is contained in $N_1^\circ$, then $W_1^\circ$ is reducible, and hence prime since it has nonempty boundary.
    Otherwise, if both $S_2$ and $S_3$ lie in $N_2^\circ$, then $W_2^\circ$ is reducible.
    Then $W_2^\circ$ must be non-prime, since it either has an $\RP^3$ summand, or contains two non-isotopic embedded spheres and hence is not $S^1 \times S^2$.
\end{proof}

\subsection{Exceptional cases of Theorem \ref{mainthm:BH_counterexamples}}\label{subsec:exceptional_cases}

We will say that the cover $p:M \rightarrow N$ is \emph{exceptional} if $N^\circ \cong W_1^\circ \# W_2^\circ$ where $W_1^\circ$ is a Seifert fibration over the disk where each branch meridian is a regular fiber and $W_2^\circ$ is another such fibration or a prime closed $3$-manifold.
These are precisely the reducible cases not covered by Theorem \ref{mainthm:BH_counterexamples}.

There is an example of an exceptional cover for which $\Ker(\L_p)$ is finite, namely the cover $p_2:S^1 \times S^2 \rightarrow S^3$ defined in Section \ref{sec:he_cover} (see Section \ref{subsec:two_factor_case} for a computation of $\Ker(\L_{p_2})$).
In this example, each $W_i^\circ$ is a solid torus.

On the other hand, we can build examples of exceptional covers for which $\Ker(\L_p)$ is infinite.
Namely, let $W_1^\circ$ be a solid torus, and let $W_2^\circ$ be any prime closed $3$-manifold where $\pi_1(W_2^\circ)$ has an element $w$ with no central power.
Then we can construct a slide homeomorphism $f$ which acts on $\pi_1(N^\circ) \cong \Z * \pi_1(W_2^\circ)$ by partially conjugating the $\pi_1(W_2^\circ)$-factor by $w$; since $w$ has no central power, the mapping class $[f] \in \Mod(N,C)$ has infinite order.
Let $p:M \rightarrow N$ be a cyclic branched cover whose monodromy homomorphism $\pi_1(N^\circ) \rightarrow \Z/k\Z$ maps $\pi_1(W_1^\circ) \cong \Z$ onto $\Z/k\Z$ and maps $\pi_1(W_2^\circ)$ to $0$.
Then $p^{-1}(W_1^\circ)$ is a $3$-sphere, and $\pi_1(M)$ is the $k$-fold free product $\pi_1(W_2^\circ) * \cdots * \pi_1(W_2^\circ)$.
Replacing $f$ with a power if necessary, we can assume that $f$ is liftable.
It has a lift $\widetilde{f}$ which conjugates each $\pi_1(W_2^\circ)$-factor by a power of $w$, and thus is overall an inner automorphism.
By Proposition \ref{prop:lifting_map_factors}, $\L_p(\alpha)$ has finite order.

In general, if a liftable homeomorphism $f$ acts on $\pi_1(N^\circ) \cong \pi_1(W_1^\circ) * \pi_1(W_2^\circ)$ by partial conjugating $\pi_1(W_i^\circ)$ by an element $w$, then its lift will act on $\pi_1^{\rel}(M) \cong H_1 * H_2$ by partially conjugating $H_i$ by a lift of $w$.
Thus, an infinite order partial conjugation of $\pi_1(N^\circ)$ will lift to a finite order automorphism of $\pi_1^{\rel}(M)$ only if $w$ lifts to an element $\widetilde{w} \in H_i$ with a central power.
However, if $W_i^\circ$ is closed or a Seifert fibration over the disk where each branch meridian is a regular fiber, then any element $w \in \pi_1(W_i^\circ)$ with no central power will lift to an element $\widetilde{w} \in H_i$ with no central power.
This is because Proposition \ref{prop:rel_pi1_is_free_product} tells us that in the first case, $H_i = \pi_1(W_i^\circ)$, and in the second case, the map $\pi_1(W_i^\circ) \rightarrow H_i$ is a quotient by a subgroup of the center $Z(\pi_1(W_i^\circ))$.
Thus, if one were to show that $\Ker(\L_p)$ is infinite for an exceptional cover $p$, then one must either study the action on $\pi_1(M)$ instead of $\pi_1^{\rel}(M)$ (as in the paragraph above), or construct a new type of element of $\Ker(\L_p)$ which is not a slide (as in Lemma \ref{lem:SF_kernel_infinite}).

\section{The hyperelliptic cover}\label{sec:he_cover}
For the remainder of this paper, we study the branched cover $p_n$ defined in the introduction, with the goal of proving Theorem \ref{mainthm:hyperelliptic_kernel} and Corollary \ref{maincor:explicit_kernel_computation}.
Here we recall the definition of $p_n$ and establish some notation.
We let $C_n \subseteq S^3$ be the $n$-component unlink, and we let $p_n:M_n \rightarrow S^3$ be the double cover branched over $C_n$.
Recall that $\pi_1(S^3 - C_n)$ is a free group on $n$ generators; the cover $p_n$ is then defined by the surjection $\varphi:\pi_1(S^3 - C_n) \rightarrow \Z/2\Z$ that maps each generator to $[1]$.
The manifold $M_n$ is the $(n-1)$-fold connected sum $(S^1 \times S^2)^{\#(n-1)}$.
The deck group $G \cong \Z/2\Z$ is generated by an involution $\tau \in \Diff^+(M_n)$.
We denote the lifting map associated to $p_n$ by 
\begin{equation*}
    \L_n:\LMod_n(S^3) \rightarrow \SMod_n(M_n)/\la \tau \ra.
\end{equation*}

In this section, we establish some basic facts about the cover $p_n$.
In Section \ref{subsec:model_of_cover}, we give a geometric description of the cover $p_n$.
In Section \ref{subsec:liftable_mcg_unlink}, we describe the liftable mapping class group $\LMod_n(S^3)$. 
We show that it is isomorphic to the group $\SymOut(F_n)$ and give an explicit presentation for this group.
In Section \ref{subsec:symmetric_mcg_unlink}, we briefly describe the symmetric mapping class group $\SMod_n(M_n)$.

\subsection{A model of the cover}\label{subsec:model_of_cover}
We construct $M_n$ as the double of $S^3$ minus $n$ $3$-balls.
That is, let $W_1$ and $W_2$ each be the complement of $n$ open $3$-balls in $S^3$.
Let $\Sigma_1^+, \ldots, \Sigma_n^+$ be the boundary components of $W_1$, and let $\Sigma_1^-, \ldots, \Sigma_n^-$ be the boundary components of $W_2$.
Then, we obtain $M_n$ by gluing each $\Sigma_i^+$ to $\Sigma_i^-$ via an orientation-reversing homeomorphism.
We let $\Sigma_i \subseteq M_n$ denote the image of $\Sigma_i^\pm$ in $M_n$.

\begin{figure}[ht]
    \centering
    \includegraphics*[scale=.35]{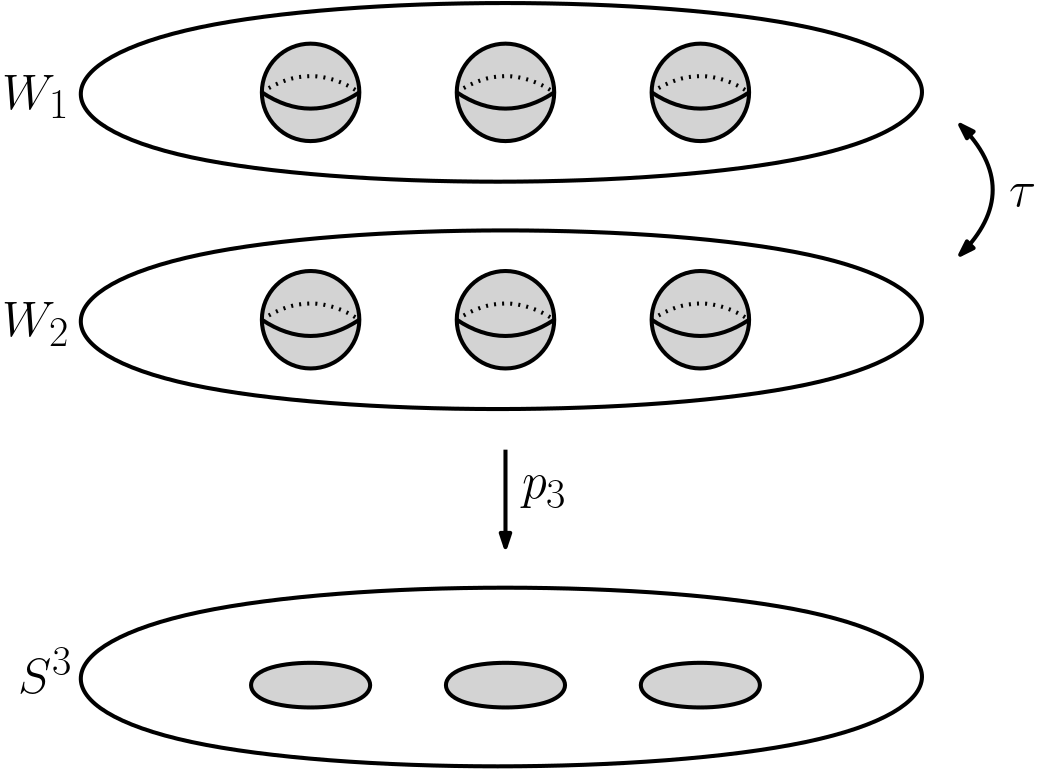}
    \caption{The cover $p_n$ in the case $n=3$.  The top hemisphere of each sphere in $W_1$ is glued to the bottom hemisphere of the opposing sphere in $W_2$.}
    \label{fig:cover_model}
\end{figure}

In order to construct $p_n$, we will directly define the involution $\tau$, and then take $p_n$ to be the quotient map.
View $S^3$ as $\R^3 \cup \{\infty\}$. 
Assume that the boundary spheres $\Sigma_i^\pm$ are Euclidean spheres centered on the $x$-axis, and that $\Sigma_i^+$ is glued to $\Sigma_i^-$ via a reflection in the $xy$-plane.
We then define $\tau$ by swapping $W_1$ and $W_2$ via $\id_{S^3}$.
The fixed set $B_n$ of $\tau$ is precisely the equator of each $\Sigma_i$ in the $xy$-plane.
The quotient $M_n/\la \tau \ra$ is homeomorphic to $S^3$, and each sphere $\Sigma_i$ descends to a disk bounded by a circle in the $xy$-plane; the union of these circles in the unlink $C_n$.
Note that $p_n$ restricts to a homeomorphism $B_n \cong C_n$.
See Figure \ref{fig:cover_model} for an illustration.

\subsection{The liftable mapping class group}\label{subsec:liftable_mcg_unlink}
Recall that $\pi_1(S^3 - C_n)$ is isomorphic to the free group $F_n$.
We fix a free generating set $y_1, \ldots, y_n$, where each $y_i$ passes once through the $i$th component of $C_n$.
Let $\SymAut(F_n)$ denote the group of \emph{symmetric automorphisms} of $F_n$, i.e.\ automorphisms that map each generator to a conjugate of another generator.
Let $\SymOut(F_n)$ denote the image of $\SymAut(F_n)$ in $\Out(F_n)$.

\begin{proposition}\label{prop:Mod_is_SymOut}
    The natural map $\Phi:\Mod(S^3,C_n) \rightarrow \Out(\pi_1(S^3-C_n))$ corestricts to an isomorphism
    \begin{equation*}
        \Mod(S^3, C_n) \cong \SymOut(F_n),
    \end{equation*}
    and the liftable mapping class group $\LMod_{p_n}(S^3)$ is the full mapping class group $\Mod(S^3, C_n)$.
\end{proposition}

Before proving Proposition \ref{prop:Mod_is_SymOut}, we describe the structure of the group $\SymOut(F_n)$.
The group $\SymOut(F_n)$ splits as a semidirect product
\begin{equation*}
    \SymOut(F_n) \cong \PSymOut(F_n) \rtimes ((\Z/2\Z)^n \rtimes S_n).
\end{equation*}
Here $\PSymOut(F_n)$ is the subgroup of \emph{pure} symmetric outer automorphisms, i.e.\ symmetric outer automorphisms that map each generator $y_i$ to a conjugate of itself.
The factor $(\Z/2\Z)^n$ corresponds to the elements $\rho_i$ that invert each generator, and the symmetric group $S_n$ corresponds to the elements that permute the generators.
By Proposition \ref{prop:Mod_is_SymOut}, this splitting corresponds to a splitting 
\begin{equation*}
    \Mod(S^3, C_n) \cong \PMod(S^3, C_n) \rtimes ((\Z/2\Z)^n \rtimes S_n).
\end{equation*}
Here $\PMod(S^3,C_n)$ is the \emph{pure mapping class group}, i.e.\ the subgroup of mapping classes that preserve the order and orientation of the components of $C_n$. The subgroup $(\Z/2\Z)^n$ reverses the orientation of the components of $C_n$, and the subgroup $S_n$ permutes the components of $C_n$.

The group $\SymAut(F_n)$ has a simple presentation, as shown by McCool \cite{mccool} as well as Gilbert \cite{gilbert} following the work of Fouxe-Rabinovitch \cite{fouxe-rabinovitch}.
From this, we also get a presentation of $\SymOut(F_n)$.
The subgroup $\PSymAut(F_n)$ is generated by the elements $\alpha_{i,j}$ for $1 \leq i \neq j \leq n$, where $\alpha_{i,j}$ maps $y_i$ to $y_jy_iy_j^{-1}$ and fixes each other generator.
This subgroup then has the defining relations
\begin{align*}
    [\alpha_{i,j},\alpha_{k,\ell}] &= 1 \\
    [\alpha_{i,k}, \alpha_{j,k}] &= 1 \\
    \alpha_{i,j}\alpha_{j,k}\alpha_{i,k} &= \alpha_{i,k}\alpha_{j,k}\alpha_{i,j},
\end{align*}
where the indices $i$, $j$, $k$, and $\ell$ are pairwise nonequal.
To complete this to a presentation of $\SymAut(F_n)$, we add order 2 generators $\rho_i$ for $1 \leq i \leq n$, where $\rho_i$ inverts $y_i$, as well as generators and relations for the symmetric group $S_n$.
The element $\rho_k$ acts on $\PSymAut(F_n)$ by 
\begin{equation*}
    \rho_k \alpha_{i,j} \rho_k^{-1} =
    \begin{cases}
        \alpha_{i,j}^{-1} & k = j \\
        \alpha_{i,j} & \text{otherwise}
    \end{cases}
\end{equation*}
and any $\sigma \in S_n$ acts on $\PSymAut(F_n)$ by 
\begin{equation*}
    \sigma \alpha_{i,j} \sigma^{-1} = \alpha_{\sigma(i),\sigma(j)}.
\end{equation*}
Finally, to obtain a presentation for $\SymOut(F_n)$, we simply quotient by the subgroup of inner automorphisms.
This amounts to adding the relation 
\begin{equation*}
    \prod_{i=1}^n \alpha_{i,j} = 1
\end{equation*}
for each $1 \leq j \leq n$.

Now, we can prove Proposition \ref{prop:Mod_is_SymOut}.
The main task is to prove that the map $\Phi:\Mod(S^3, C_n) \rightarrow \Out(F_n)$ is an isomorphism onto $\SymOut(F_n)$.
The fact that $\LMod_{p_n}(S^3) = \Mod(S^3, C_n)$ then follows directly from the fact that any symmetric automorphism of $\pi_1(S^3 - C_n) \cong F_n$ preserves the monodromy homomorphism $\varphi:\pi_1(S^3 - C_n) \rightarrow \Z/2\Z$ of the cover $p_n$.

The isomorphism follows from results about the \emph{(extended) loop braid group}.
This is the group 
\begin{equation*}
    \Mod(B^3, C_n) \coloneqq \pi_0(\Homeo^+_{\del}(B^3, C_n)),
\end{equation*}
where $\Homeo^+_{\del}(B^3, C_n)$ is the group of homeomorphisms of the pair $(B^3, C_n)$ which fix $\del B^3$ pointwise and preserve the orientation of $B^3$ (but not necessarily the orientation of $C_n$).
The action on $\pi_1$ induces a natural map $\widehat{\Phi}:\Mod(B^3, C_n) \rightarrow \Aut(F_n)$.
Goldsmith \cite{goldsmith}, following the work of Dahm \cite{dahm}, proved that this map corestricts to an isomorphism
\begin{equation*}
    \Mod(B^3, C_n) \cong \SymAut(F_n)
\end{equation*}
(Goldsmith worked with a different definition of the loop braid group in terms of \emph{motion groups}, see Damiani's survey \cite{damiani} for a variety of definitions and their relations).

From this isomorphism, we can deduce Proposition \ref{prop:Mod_is_SymOut}.
\begin{proof}[Proof of Proposition \ref{prop:Mod_is_SymOut}]
    Let $S = \del B^3$, and view $S^3$ as being obtained from $B^3$ by capping $S$ with a $3$-ball $E$.
    Observe that we have a map $\F:\Mod(B^3, C_n) \rightarrow \Mod(S^3, C_n)$ by extending any homeomorphism of $(B^3, C_n)$ to be the identity on $E$.
    
    First, we claim that the map $\F$ is surjective.
    Take any $\alpha \in \Mod(S^3, C_n)$, and choose a representative $f$ which preserves the $3$-ball $E$.
    Since $f$ preserves the orientation of $S^3$, $f$ must preserve the orientation of $E$.
    Since $\Mod(S^2)$ is trivial, we may assume up to isotopy that $f$ fixes $S$ pointwise.
    Then by the Alexander trick, we may assume up to isotopy that $f$ fixes $E$ pointwise.
    Thus we have $[f] = \F([f\v_{B^3}])$.

    Now, we can prove the proposition.
    We have the following commutative diagram:
    \begin{equation*}
        \begin{tikzcd}
            {\Mod(B^3,C_n)} & {\SymAut(F_n)} \\
            {\Mod(S^3,C_n)} & {\Out(F_n)}
            \arrow["{\widehat{\Phi}}"', "\cong", from=1-1, to=1-2]
            \arrow["\F"', two heads, from=1-1, to=2-1]
            \arrow[from=1-2, to=2-2]
            \arrow["\Phi"', from=2-1, to=2-2]
        \end{tikzcd}
    \end{equation*}
    This diagram immediately implies that $\Im(\Phi) = \SymOut(F_n)$.
    To see that $\Phi$ is injective, suppose $\alpha \in \Mod(S^3, C_n)$ and $\Phi(\alpha)$ is trivial.
    Then $\alpha = \F(\beta)$ where $\beta \in \Mod(B^3, C_n)$ maps to an inner automorphism under $\widehat{\Phi}$.
    Note that any inner automorphism of $F_n$ is realized by a homeomorphism of $(B^3, C_n)$ obtained by sliding $S$ along a loop based at $* \in S$ (see Section \ref{subsec:slide_homeos} for a precise definition of a slide homeomorphism).
    Since $\widehat{\Phi}$ is injective, we know that $\beta$ must be represented by such a slide.
    But such a slide extends to a trivial element of $\Mod(S^3,C_n)$, and so $\alpha$ must be trivial.
\end{proof}

\subsection{The symmetric mapping class group}\label{subsec:symmetric_mcg_unlink}
Next, we briefly comment on the group $\SMod_n(M_n)$.
Laudenbach \cite{laudenbach-1,laudenbach-2} proved that the action of $\Mod(M_n)$ on $\pi_1(M_n) \cong F_{n-1}$ yields a short exact sequence
\begin{equation*}
    1 \rightarrow \Twist(M_n) \rightarrow \Mod(M_n) \rightarrow \Out(F_{n-1}) \rightarrow 1,
\end{equation*}
where $\Twist(M_n)$ is the subgroup generated by all sphere twists.
Laudenbach also showed that $\Twist(M_n) \cong (\Z/2\Z)^{n-1}$ and is generated by twists around the core spheres $\Sigma_1, \ldots, \Sigma_{n-1}$ in our model of $M_n$.
We note that Brendle-Broaddus-Putman \cite{brendle-broaddus-putman} proved that this sequence splits, though we will not use this fact.
The subgroup $\SMod_n(M_n) \leq \Mod(M_n)$ lies in the normalizer of $\la \tau \ra$; this normalizer is in fact the centralizer since $\tau$ has order $2$.
Moreover, since $\tau$ preserves each $\Sigma_i$, the subgroup $\Twist(M_n)$ also lies in the centralizer of $\tau$.
Thus we have an exact sequence 
\begin{equation*}
    1 \rightarrow \Twist(M_n) \rightarrow \SMod_n(M_n) \rightarrow \Out_{\tau_*}(F_{n-1}),
\end{equation*}
where $\tau_* \in \Out(F_{n-1})$ is the outer automorphism induced by $\tau$, and $\Out_{\tau_*}(F_{n-1})$ is the centralizer of $\tau_*$.
The surjectivity of $\Mod_{n}(M_n) \rightarrow \Out(F_{n-1})$ does not immediately imply that the map $\SMod_n(M_n) \rightarrow \Out_{\tau_*}(F_{n-1})$ is surjective; it would be of interest to determine whether this map is indeed surjective.

\begin{remark}\label{rem:palindromic_automorphisms}
    As pointed out by Margalit-Winarski \cite[\S 11]{margalit-winarski}, there is a connection between $\SMod_n(M_n)$ and the group $\PalAut(F_{n-1})$ of \emph{palindromic automorphisms} of $F_{n-1}$.
    These are the automorphisms of $F_{n-1}$ which map each generator $x_i$ to a word of the form $wx_j^{\pm 1}\overline{w}$, where $x_j$ is another generator and $\overline{w}$ is the reverse of $w$.
    The group $\PalAut(F_{n-1})$ has a presentation similar to that of $\SymAut(F_n)$, as shown by Collins \cite{collins}.

    The connection is as follows.
    Let $\iota \in \Aut(F_{n-1})$ be the automorphism that inverts each generator.
    Then $\PalAut(F_{n-1})$ is precisely the centralizer of $\iota$ in $\Aut(F_{n-1})$.
    Moreover, the outer automorphism $[\iota] \in \Out(F_{n-1})$ is precisely the outer automorphism $\tau_*$.
    Thus the projection $\Aut(F_{n-1}) \rightarrow \Out(F_{n-1})$ restricts to a map $\PalAut(F_{n-1}) \rightarrow \Out_{\tau_*}(F_{n-1})$.
    This map is in fact injective since no inner automorphism of $F_{n-1}$ lies in $\PalAut(F_{n-1})$.
    However, we caution the reader this map is \emph{not} surjective (for $n \geq 3$).
    For example, let $x_1, \ldots, x_{n-1}$ be the generators of $F_{n-1}$, and let $\sigma \in \Aut(F_{n-1})$ be the automorphism which maps $x_1$ to $x_1^{-1}$ and each other generator $x_i$ to $x_ix_1^{-1}$.
    A direct calculation shows that the commutator $[\sigma, \iota]$ is inner, and so $\sigma$ projects to an element of $\Out_{\tau_*}(F_{n-1})$.
    Howover, any representative of the outer automorphism $[\sigma]$ maps each $x_i$ to an even length word for $i \neq 1$.
    Thus $[\sigma]$ is not represented by any palindromic automorphism.
\end{remark}

\section{Algebraic characterization of the lifting kernel}\label{sec:alg_characterization}
Our ultimate goal is to prove Theorem \ref{mainthm:hyperelliptic_kernel}, i.e.\ to compute the kernel of the lifting map
\begin{equation*}
    \L_n:\Mod(S^3,C_n) \rightarrow \SMod_n(M_n)/\la \tau_* \ra.
\end{equation*}
In this section, we will recast this problem into a purely algebraic one by using the action of $\Mod(S^3,C_n)$ on the relative fundamental group $\pi_1^{\rel}(M_n)$.

First, note that by post-composing $\L_n$ with a quotient map, we get a map
\begin{equation*}
    \overline{\L}_n:\Mod(S^3, C_n) \rightarrow \frac{\SMod_n(M_n)/\la \tau_* \ra}{\Twist(M_n)}.
\end{equation*}
Next, let $H_n$ denote the $n$-fold free product $(\Z/2\Z)^{*n}$.
Then we have a natural projection $\pi:F_n \rightarrow H_n$ by reducing each generator mod 2.
We define $\SymAut(H_n)$ and $\SymOut(H_n)$ just as with $F_n$.
Since any symmetric automorphism of $F_n$ preserves $\Ker(\pi)$, we get a map
\begin{equation*}
    \P_n:\SymOut(F_n) \rightarrow \SymOut(H_n).
\end{equation*}
The main result of this section gives a correspondence between $\Ker(\overline{\L}_n)$ and $\Ker(\P_n)$.
Recall from Proposition \ref{prop:Mod_is_SymOut} that the action of $\Mod(S^3,C_n)$ on $\pi_1(S^3 - C_n)$ induces an isomorphism $\Mod(S^3,C_n) \cong \SymOut(F_n)$.

\begin{proposition}\label{prop:lifting_kernel_is_alg_kernel}
    For $n \geq 3$, the isomorphism $\Mod(S^3, C_n) \cong \SymOut(F_n)$ maps $\Ker(\overline{\L}_n)$ isomorphically onto $\Ker(\P_n)$.
\end{proposition}

\noindent We note that Proposition \ref{prop:lifting_kernel_is_alg_kernel} is false in the case $n=2$, since $\Ker(\overline{\L}_2)$ is all of $\Mod(S^3,C_2)$; see Section \ref{subsec:two_factor_case} for details.

We will prove Proposition \ref{prop:lifting_kernel_is_alg_kernel} in three steps.
First, in Section \ref{subsec:rel_pi1_he_cover}, we show that $\pi_1^{\rel}(M_n)$ is isomorphic to $H_n$ and identify the subgroup $\pi_1(M_n) \leq \pi_1^{\rel}(M_n)$.
In Section \ref{subsec:restriction_map}, we show that the natural restriction map 
\begin{equation*}
    \overline{r}:\Out(\pi_1^{\rel}(M_n),\pi_1(M)) \rightarrow \Out_{\tau_*}(\pi_1(M))/\la \tau_* \ra
\end{equation*}
is injective (the map $\overline{r}$ was defined in Section \ref{subsec:factoring_lifting_map}).
Finally, in Section \ref{subsec:proof_lifting_is_alg_kernel}, we complete the proof of Proposition \ref{prop:lifting_kernel_is_alg_kernel} using the factorization of $\L_n$ from Proposition \ref{prop:lifting_map_factors}.

\subsection{The relative fundamental group}\label{subsec:rel_pi1_he_cover}
We begin by identifying the relative fundamental group $\pi_1^{\rel}(M_n)$ and its subgroup $\pi_1(M_n)$.
Let $\varphi:\pi_1(S^3 - C_n) \rightarrow G$ be the monodromy homomorphism for the cover $p_n$.
This is equivalent the map $\varphi':F_n \rightarrow \Z/2\Z$ defined by mapping each generator to $[1]$.
We define a map $\psi:H_n \rightarrow \Z/2\Z$ by mapping each generator to $[1]$.
Then we have the following.

\begin{lemma}\label{lem:hyperelliptic_rel_htpy_groups}
    We have the following commutative diagram:
    \begin{equation*}
        \begin{tikzcd}
            & {} & {\pi_1(M_n^\circ)} & {\pi_1(S^3 - C_n)} & G \\
            {F_{2n-1}} & {F_n} & {\Z/2\Z} && {\pi_1(M_n)} & {\pi_1^{\rel}(M_n)} & {G} \\
            && {F_{n-1}} & {H_n} & {\Z/2\Z}
            \arrow[hook, from=1-3, to=1-4]
            \arrow[two heads, "\varphi", from=1-4, to=1-5]
            \arrow[shorten <=12pt, shorten >=24pt, Leftrightarrow, from=1-4, to=2-2]
            \arrow[shorten <=18pt, shorten >=18pt, Rightarrow, from=1-4, to=2-6]
            \arrow[hook, from=2-1, to=2-2]
            \arrow[two heads, "\varphi'"{pos=0.4}, from=2-2, to=2-3]
            \arrow[shorten <=18pt, shorten >=18pt, Rightarrow, from=2-2, to=3-4]
            \arrow[hook, from=2-5, to=2-6]
            \arrow[two heads, from=2-6, to=2-7]
            \arrow[shorten <=18pt, shorten >=18pt, Leftrightarrow, from=2-6, to=3-4]
            \arrow[hook, from=3-3, to=3-4]
            \arrow[two heads, "\psi"', from=3-4, to=3-5]
        \end{tikzcd}
    \end{equation*}
    Here each sequence of $3$ groups is short exact.
    The arrow $\Rightarrow$ denotes a surjection of short exact sequences, while the arrow $\Leftrightarrow$ denotes an isomorphism of short exact sequences.
\end{lemma}
\begin{proof}
    First, we can verify the existence of the four short exact sequences.
    \begin{itemize}
        \item The top and right short exact sequences are precisely the ones found in the commutative diagram (\ref{eqn:htpy_group_diagram}), which relates $\pi_1(N^\circ)$ to $\pi_1^{\rel}(M)$ for a general cover $M \rightarrow N$.
        
        \item For the left short exact sequence, we must verify that $\Ker(\varphi') \cong F_{2n-1}$.
        If we let $R_n$ denote the $n$-petaled rose, then $\varphi'$ corresponds to a double cover $\widetilde{R}_n \rightarrow R_n$ with $\pi_1(\widetilde{R}_n) \cong \Ker(\varphi')$, and
        \begin{equation*}
            \text{rank}(\pi_1(\widetilde{R}_n))
            = 1 - \chi(\widetilde{R}_n)
            = 1 - 2\chi({R}_n)
            = 1 - 2(1 - n)
            = 2n - 1.
        \end{equation*}

        \item For the bottom short exact sequence, we must verify that $\Ker(\psi) \cong F_{n-1}$.
        To see this, let $z_1, \ldots, z_n$ be a basis of $H_n$, and let $T$ be the Cayley graph with respect to this generating set.
        Then $\Ker(\psi)$ acts freely on the vertices of $T$.
        Note that $\Ker(\psi)$ is precisely the set of even length words in the letters $z_i$.
        It follows that $\Ker(\psi)$ has precisely two vertex orbits on $T$, corresponding to the words of even and odd length.
        Moreover, for each edge of $T$, the endpoints lie in different orbits, and so $\Ker(\psi)$ acts freely on the set of edges of $T$.
        It follows that $T/\Ker(\psi)$ is a graph with two vertices and $n$ edges between them, and thus $\Ker(\psi)$ is free of rank $n-1$.
    \end{itemize}

    Next, we can identify the maps between the short exact sequences.  The commutativity of the diagram will follow directly.

    \begin{itemize}
        \item The top-right surjection is precisely the surjection of short exact sequences of homotopy groups in the commutative diagram (\ref{eqn:htpy_group_diagram}).
        \item The top-left isomorphism is induced by the isomorphism $\pi_1(S^3 - C_n) \cong F_n$.
        \item The bottom-left surjection is induced by natural surjection $\pi:F_n \rightarrow H_n$.
        \item For the bottom-right isomorphism, we claim there is an isomorphism $\pi_1^{\rel}(M_n) \rightarrow H_n$ taking the natural map $\pi_1^{\rel}(M_n) \rightarrow G$ to the map $\psi$.  This follows from Proposition \ref{prop:rel_pi1_is_free_product}.  Indeed, the branch meridians of $S^3 - C_n$ are precisely the free generators of $\pi_1(S^3 - C_n)$, and they each map to the order $2$ generator of $G$, and thus $\pi_1^{\rel}(M_n)$ is the quotient of $\pi_1(S^3 - C_n) \cong F_n$ by the normal closure of the squares of the generators.
    \end{itemize}
\end{proof}

\subsection{The restriction map}\label{subsec:restriction_map}
Next, as in Section \ref{subsec:factoring_lifting_map}, let
\begin{equation*}
    \overline{r}:\Out(\pi_1^{\rel}(M_n),\pi_1(M_n)) \rightarrow \Out_{\tau_*}(\pi_1(M_n))/\la \tau_* \ra
\end{equation*}
be the map induced from the restriction map $\Aut(\pi_1^{\rel}(M_n),\pi_1(M_n)) \rightarrow \Aut(\pi_1(M_n))$.
By Lemma \ref{lem:hyperelliptic_rel_htpy_groups}, the map $\overline{r}$ can be equivalently seen as a map 
\begin{equation*}
    \Out(H_n, F_{n-1}) \rightarrow \Out_{\tau_*}(F_{n-1})/\la \tau_* \ra,
\end{equation*}
where $F_{n-1}$ is the kernel of the map $\psi:H_n \rightarrow \Z/2\Z$ defined in Section \ref{subsec:rel_pi1_he_cover}.
Then, we can restrict this to a map 
\begin{equation*}
    s:\SymOut(H_n) \rightarrow \Out_{\tau_*}(F_{n-1})/\la \tau_* \ra.
\end{equation*}
Our next lemma shows that this map is injective.

\begin{lemma}\label{lem:restriction_map}
    For $n \geq 3$, the restriction map $s:\SymOut(H_n) \rightarrow \Out_{\tau_*}(F_{n-1})/ \la \tau_* \ra$ is injective.
\end{lemma}

\noindent Lemma \ref{lem:restriction_map} is false in the case $n=2$, since the outer automorphism that swaps the two generators of $H_n$ will descend to $\tau_*$.

\begin{proof}[Proof of Lemma \ref{lem:restriction_map}]
    Let $z_1, \ldots, z_n$ be a basis of $H_n$.
    We show first that the elements $x_i \coloneqq z_iz_n$ for $1 \leq i < n$ generate the subgroup $\Ker(\psi) \cong F_{n-1}$.
    This is because $\Ker(\psi)$ is precisely the set of even length words in the letters $z_i$, and for any even length word $w = z_{i_1}z_{j_1} \cdots z_{i_m}z_{j_m}$, we can replace each $z_{i_k}z_{j_k}$ with the word $z_{i_k}z_nz_nz_{j_k}$ to get a word on the elements $x_i = z_iz_n$.

    Just as with the group $\SymOut(F_n)$, we have a splitting $\SymOut(H_n) \cong \PSymOut(H_n) \rtimes S_n$, where $\PSymOut(H_n)$ is the subgroup of outer automorphisms sending each generator to a conjugate of itself, and where the symmetric group permutes the generators (unlike $\SymOut(F_n)$, we don't get a $(\Z/2\Z)^n$-factor since each generator has order $2$).
    We will show that $s$ is injective on $\PSymOut(H_n)$ and on $S_n$, and then check that the images of these subgroups under $r$ have trivial intersection.

    Suppose first that $\alpha \in \PSymOut(H_n)$, and that $s(\alpha)$ is trivial.
    We can pick a representative $f$ of $\alpha$ such that $f(z_{n}) = z_{n}$.
    Then for each $i < n$, we have that $f(z_i) = w_iz_iw_i^{-1}$ for some $w_i \in H_n$.
    Note that either $w_i \in F_{n-1}$ or $w_iz_n \in F_{n-1}$.
    We can compute that
    \begin{equation*}
        f(x_i) =
        \begin{cases}
            w_ix_i\overline{w}_i & \text{if } w_i \in F_{n-1} \\
            w_iz_{n}x_i^{-1}\overline{z_nw_i} & \text{if } w_iz_n \in F_{n-1}
        \end{cases}
    \end{equation*}
    where $\overline{z_nw_i}$ is the reverse of $z_nw_i$ with respect to the generating set $x_1, \ldots, x_{n-1}$ of $F_{n-1}$.
    Recall that $\tau_* = \iota_*$, where $\iota \in \Aut(F_{n-1})$ inverts each $x_i$.
    Since $s(\alpha)$ is trivial, we know that either $f$ or $f\iota$ is an inner automorphism of $F_{n-1}$.
    This is possible only if all $w_i$ are trivial, or each $w_i = z_{n}$.
    But in either case, this implies that $\alpha$ is trivial.

    Next, suppose that $\sigma \in S_{n}$.
    Note that if $\sigma(z_n) = z_n$, then $\sigma$ simply permutes the indices of the generators $x_1, \ldots, x_{n-1}$.
    Otherwise, we can write $\sigma = (i\ n) \circ \sigma'$ for some $1 \leq i < n$, where $(i\ n)$ is the transposition swapping $z_i$ and $z_n$, and $\sigma'$ is some permutation that fixes $z_n$.
    The element $(i\ n)$ maps $x_i$ to $x_i^{-1}$ and each other $x_j$ to $x_jx_i^{-1}$, and $\sigma'$ permutes the indices of the generators $x_1, \ldots, x_{n-1}$.
    In any case, we see that if $\sigma$ is nontrivial, then $s(\sigma)$ is nontrivial.

    Finally, it follows from our computations that $s(\alpha) \neq s(\sigma)$ for any nontrivial $\alpha \in \PSymOut(H_n)$ and $\sigma \in S_{n}$.
\end{proof}

\begin{remark}\label{rem:rel_pi1_and_braids_discussion}
    Lemmas \ref{lem:hyperelliptic_rel_htpy_groups} and \ref{lem:restriction_map} provide new insight into the connection between the lifting map and Magnus' question as discussed in Section \ref{subsec:rel_pi1_and_braid}.
    Recall that the question was whether the map 
    \begin{equation*}
        \eta_{n,k}:B_n \rightarrow \SymAut(H_{n,k})
    \end{equation*}
    is injective, where $H_{n,k}$ is $n$-fold free product $(\Z/k\Z)^{*n}$.
    
    Birman-Hilden \cite{birman-hilden} proved that the map $\eta_{n,k}$ is injective; their argument is as follows (see also \cite[\S 7]{margalit-winarski}).
    Let $D_n$ denote the $n$-times punctured disk.
    Recall that $\pi_1(D_n)$ is the free group $F_n$, and is generated by the elements $x_1, \ldots, x_n$ where $x_i$ is a loop around the $i$th puncture.
    Let $p_{n,k}^\circ:S_{n,k}^\circ \rightarrow D_n$ denote the cover corresponding to the map $\varphi_k:F_n \rightarrow \Z/k\Z$ that maps each $x_i$ to $[1]$; this completes to a branched cover of the disk $p_{n,k}:S_{n,k} \rightarrow D$.
    Then $x_i^k$ lies in $\pi_1(S_{n,k}^\circ) \leq \pi_1(D_n)$ for each $i$, and Van Kampen's theorem implies that the kernel of $\pi_1(S_{n,k}^\circ) \rightarrow \pi_1(S)$ is the normal closure $N_{n,k}$ of the set $\{x_1^k, \ldots, x_n^k\}$.
    Since $F_n/N_{n,k} \cong H_{n,k}$, we get the following commutative diagram with exact rows:
    \begin{equation}\label{eqn:braid_group_diagram}
        \begin{tikzcd}
            1 & {\pi_1(S^\circ)} & {F_n} & {\Z/k\Z} & 1 \\
            1 & {\pi_1(S)} & {H_{n,k}} & {\Z/k\Z} & 1
            \arrow[from=1-1, to=1-2]
            \arrow[from=1-2, to=1-3]
            \arrow[two heads, from=1-2, to=2-2]
            \arrow["{\varphi_k}", from=1-3, to=1-4]
            \arrow[two heads, from=1-3, to=2-3]
            \arrow[from=1-4, to=1-5]
            \arrow[equals, from=1-4, to=2-4]
            \arrow[from=2-1, to=2-2]
            \arrow[from=2-2, to=2-3]
            \arrow[from=2-3, to=2-4]
            \arrow[from=2-4, to=2-5]
        \end{tikzcd}
    \end{equation}
    We have a lifting map $\L_{n,k}:B_n \rightarrow \SMod(S_{n,k})$ (we don't quotient by the deck group since we can fix a base point on the boundary).
    Birman-Hilden's work implies that $\L_{n,k}$ is injective, which means that $B_n$ acts faithfully on the subgroup $\pi_1(S_{n,k}) \leq H_{n,k}$, and hence it acts faithfully on $H_{n,k}$ as well.

    As noted in Section \ref{subsec:rel_pi1_and_braid}, the application of Birman-Hilden's result is initially surprising, but made more natural using the relative fundamental group.
    Namely, following the proofs of Proposition \ref{prop:rel_pi1_is_free_product} and Lemma \ref{lem:hyperelliptic_rel_htpy_groups}, the relative fundamental group $\pi_1^{\rel}(S_{n,k})$ is precisely the group $H_{n,k}$. 
    Furthermore, the diagram (\ref{eqn:braid_group_diagram}) is corresponds to the first and third rows of the commutative diagram (\ref{eqn:htpy_group_diagram}).
    In the case $k = 2$, we have that $\pi_1(D_n) \cong \pi_1(S^3 - C_n)$ and $\pi_1^{\rel}(S_{n,2}) \cong \pi_1^{\rel}(M_n)$, and the diagram (\ref{eqn:braid_group_diagram}) is precisely the diagram (\ref{eqn:htpy_group_diagram}) for the cover $p_n:M_n \rightarrow S^3$.

    In fact, by adapting the proof of Lemma \ref{lem:restriction_map}, one can show that in the case $k=2$, the natural map $\SymAut(H_{n,2}) \rightarrow \Aut(\pi_1(S_{n,2}))$ is injective.
    Thus, Birman-Hilden's injectivity result for $\L_{n,2}$ is in fact equivalent to the injectivity of the map $\eta_{n,2}:B_n \rightarrow \SymAut(H_{n,2})$.
    We note that there are purely algebraic proofs that $\eta_{n,k}$ is injective, see e.g.\ \cite{bacardit-dicks} or \cite{johnson}.
    It would be of interest to prove that $\eta_{n,k}$ is injective using the complexes $K_0(F_n)$ and $K_0(H_{k,n})$, as in our proof of Theorem \ref{mainthm:hyperelliptic_kernel} below.
\end{remark}

\subsection{Proof of Proposition \ref{prop:lifting_kernel_is_alg_kernel}}\label{subsec:proof_lifting_is_alg_kernel}

Now, we can prove Proposition \ref{prop:lifting_kernel_is_alg_kernel}.
By Proposition \ref{prop:lifting_map_factors}, we have a commuting square 
\begin{equation*}
    \begin{tikzcd}
        {\Mod(S^3,C_n)} & {\Out(\pi_1^{\rel}(M_n),\pi_1(M_n))} \\
        {\frac{\SMod_n(M_n)/\la \tau \ra}{\Twist(M_n)}} & {\Out_{\tau_*}(\pi_1(M_n))/\la \tau_* \ra}
        \arrow["\Psi", from=1-1, to=1-2]
        \arrow["{\overline{\L}_n}"', from=1-1, to=2-1]
        \arrow["{\overline{r}}", from=1-2, to=2-2]
        \arrow["{\overline{\Theta}}"', hook, from=2-1, to=2-2]
    \end{tikzcd}
\end{equation*}
where $\overline{\Theta}$ is the embedding induced by $\Theta:\SMod_n(M_n)/\la \tau_* \ra \rightarrow \Out_{\tau_*}(\pi_1(M_n))/\la \tau_* \ra$.
By Lemma \ref{lem:hyperelliptic_rel_htpy_groups}, this diagram is equivalent to the following diagram:
\begin{equation*}
    \begin{tikzcd}
        {\Mod(S^3,C_n)} & {\SymOut(H_n)} \\
        {\frac{\SMod_n(M_n)/\la \tau \ra}{\Twist(M_n)}} & {\Out_{\tau_*}(F_{n-1})/\la \tau_* \ra}
        \arrow["\Psi", from=1-1, to=1-2]
        \arrow["{\overline{\L}_n}"', from=1-1, to=2-1]
        \arrow["s", from=1-2, to=2-2]
        \arrow["{\overline{\Theta}}"', hook, from=2-1, to=2-2]
    \end{tikzcd}
\end{equation*}
By Lemma \ref{lem:restriction_map}, the map $s$ is injective, and hence $\Ker(\overline{\L}_n) = \Ker(\Psi)$.
But Lemma \ref{lem:hyperelliptic_rel_htpy_groups} tells us that the isomorphism $\Mod(S^3,C_n) \cong \SymOut(F_n)$ takes the map $\Psi$ to the map $\P_n$, and thus it maps $\Ker(\overline{\L}_n)$ isomorphically onto $\Ker(\P_n)$.

\section{Computing the algebraic kernel}\label{sec:alg_kernel}
Proposition \ref{prop:lifting_kernel_is_alg_kernel} tells us that the main task in proving Theorem \ref{mainthm:hyperelliptic_kernel} is to understand the kernel of the map 
\begin{equation*}
    \P_n:\SymOut(F_n) \rightarrow \SymOut(H_n),
\end{equation*}
where $H_n$ is the $n$-fold free product $(\Z/2\Z)^{*n}$.
In this section, we accomplish this by finding a normal generating set for $\Ker(\P_n)$.
Namely, for $1 \leq i \leq n$, let $\widehat{\rho}_i \in \SymAut(F_n)$ be the automorphism that inverts the $i$th generator and fixes the remaining generator, and let $\rho_i \in \SymOut(F_n)$ be the outer automorphism represented by $\widehat{\rho}_i$.
Then we have the following.

\begin{proposition}\label{prop:alg_kernel_computation}
    The kernel $\Ker(\P_n)$ is the $\SymOut(F_n)$-normal closure of the set $\{\rho_1, \ldots, \rho_n\}$.
\end{proposition}

To prove Proposition \ref{prop:alg_kernel_computation}, we employ certain simplicial complexes $K_0(F_n)$ and $K_0(H_n)$ defined by McCullough-Miller \cite{mccullough-miller}.
The vertices of these complexes are \emph{labelled bipartite trees}, and the simplices correspond to chains in a partial ordering on the set of labelled bipartite trees given by \emph{folding}.
The key property proved by McCullough-Miller is that each of these complexes are contractible.
We then show that $K_0(F_n)/\Ker(\P_n) \cong K_0(H_n)$; since both complexes are contractible, this allows us to read off generators of $\Ker(\P_n)$ from its vertex stabilizers.

In Section \ref{subsec:the_complexes}, we recall the definition of the complexes $K_0(F_n)$ and $K_0(H_n)$, as well as some basic properties.
In Section \ref{subsec:complexes_quotient_map}, we show that $K_0(F_n)/\Ker(\P_n) \cong K_0(H_n)$.
In Section \ref{subsec:pf_alg_kernel_complutation}, we complete the proof of Proposition \ref{prop:alg_kernel_computation} by computing the vertex stabilizers of $\Ker(\P_n)$.

\subsection{The simplicial complexes}\label{subsec:the_complexes}

We begin by defining the simplicial complexes, following the work of McCullough-Miller \cite{mccullough-miller}.
We note that McCullough-Miller primarily work with a definition in terms of ``based partitions''; however, we will use their alternative definition in terms of ``labelled bipartite trees'' \cite[\S 1.3]{mccullough-miller}, as this definition is similar to more familiar complexes like the Culler-Vogtman Outer space.
We also note that McCullough-Miller discuss arbitrary free products $G_1 * \cdots * G_n$, where $G_i$ is any freely indecomposable group; we are only concerned with the case that $G_i = \Z$ or $G_i = \Z/2\Z$, which will slightly simplify the discussion.

Let $L = F_n$ or $L = H_n$.
Fix a preferred basis $B_0$ of $L$.
McCullough-Miller \cite{mccullough-miller} construct a contractible simplicial complex on which $\SymOut(L)$ acts as follows.

\begin{figure}[ht]
    \centering
    \includegraphics*[scale=.30]{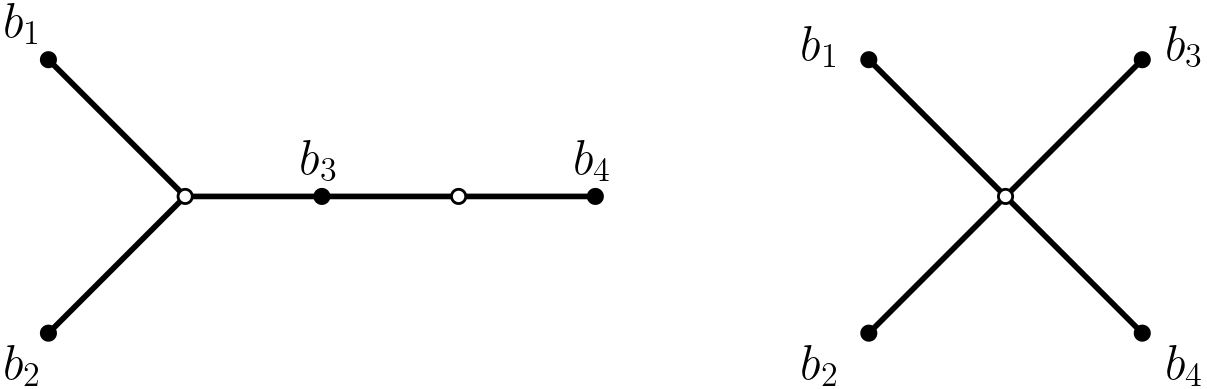}
    \caption{Two bipartite trees labelled by the basis $\{b_1, b_2, b_3, b_4\}$.}
    \label{fig:trees}
\end{figure}

We define a \emph{labelled bipartite tree} to be a triple $\T = (T, B, \lambda)$ where
\begin{itemize}
    \item $T$ is a tree,
    \item $B = \{b_1, \ldots, b_n\}$ is a basis of $L$,
    \item $\lambda:B \hookrightarrow V(T)$ is an embedding of $B$ into the vertex set of $T$,
\end{itemize}
subject to the conditions
\begin{itemize}
    \item each edge of $T$ has exactly one vertex in $\Im(\lambda)$,
    \item each vertex with valence 1 lies in $\Im(\lambda)$.
\end{itemize}
See Figure \ref{fig:trees} for examples.
We call the vertices in $\Im(\lambda)$ the \emph{labelled vertices} of $\T$.
Two labelled bipartite trees $(T,B,\lambda)$ and $(T',B',\lambda')$ are \emph{isomorphic} if $B = B'$ and there is an isomorphism of trees $T \cong T'$ taking $\lambda$ to $\lambda'$.
If $\T$ has a single unlabelled vertex, we call it a \emph{trivial $B$-labelled bipartite tree}.
In Figure \ref{fig:trees}, the graph on the right is a trivial $\{b_1,b_2,b_3,b_4\}$-labelled bipartite tree.
For each basis $B$, there is a unique trivial $B$-labelled bipartite tree up to isomorphism.
The \emph{type} of $\T$ is its isomorphism class as an (unlabelled) bipartite tree; that is, $(T,B,\lambda)$ and $(T',B',\lambda')$ have the same type if there is an isomorphism of trees $T \cong T'$ taking $\Im(\lambda)$ to $\Im(\lambda')$.

The set of labelled bipartite trees comes with a partial ordering.
If $v$ is a labelled vertex of $\T = (T, B, \lambda)$ with two adjacent edges $e_1$ and $e_2$, the \emph{fold} of $\T$ along $e_1$ and $e_2$ is $(T', B, \lambda')$, where $T'$ is the quotient of $T$ where we identify $e_1$ and $e_2$ and their unlabelled vertices, and $\lambda'$ is the labelling induced by $\lambda$.
In Figure \ref{fig:trees}, the graph on the right is obtained from the graph on the left by folding at the vertex labelled by $b_3$.
We say $\T' \leq \T$ if $\T'$ can be obtained from $\T$ by a sequence of folds.
The trivial labelled bipartite trees are the minimal elements of this poset.
Moreover, one can show that this poset has height $n-1$, and in particular a labelled bipartite tree has at most $n-1$ unlabelled vertices \cite[Lem~1.3]{mccullough-miller}.

\begin{figure}[ht]
    \centering
    \includegraphics*[scale=.35]{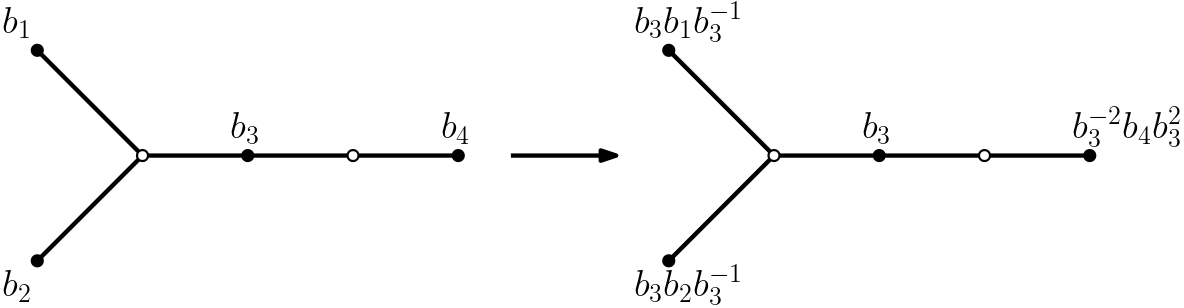}
    \caption{A vertex automorphism based at the vertex labelled by $b_3$.}
    \label{fig:vertex_aut}
\end{figure}

Observe that $\Aut(L)$ acts on the set $\B$ of bases of $L$, and hence $\Aut(L)$ acts on the set of labelled bipartite trees.
Moreover, we can see that $\Aut(L)$ preserves the partial ordering.
Given a labelled bipartite tree $\T$, a \emph{vertex automorphism carried by $\T$} is an automorphism $f_{v,\theta}$ where
\begin{itemize}
    \item $v$ is a labelled vertex, say with label $b_v \in B$,
    \item $\theta:B \rightarrow \la b_v \ra$ is a function such that $\theta(b_1) = \theta(b_2)$ if $\lambda(b_1)$ and $\lambda(b_2)$ lie on the same component of $T - v$,
\end{itemize}
whose action on $L$ is given in the basis $B$ by
\begin{equation*}
    f_{v,\theta}(b) = \theta(b)b\theta(b)^{-1}.
\end{equation*}
In other words, $f_{v,\theta}$ conjugates each $b \in B$ by some power of $b_v$, where basis elements on the same component of $T - v$ must be conjugated by the same power.
See Figure \ref{fig:vertex_aut} for an example.
We note that McCullough-Miller call $f_{v,\theta}$ a \emph{symmetric Whitehead automorphism}.

Finally, we can define the simplicial complex.
We let $\sim$ be the equivalence relation generated by the relation that $\T \sim \T'$ whenever $\T' = f_{v,\theta}\T$ for some vertex automorphism $f_{v,\theta}$ carried by $\T$.
One can check that $\sim$ descends to an equivalence relation between isomorphism classes of labelled bipartite trees.
We define 
\begin{equation*}
    \K(L) \coloneqq \{\text{isomorphism classes of labelled bipartite trees}\}/\sim.
\end{equation*}
Since $\Aut(L)$ preserves the partial ordering, the set $\K(L)$ inherits a poset structure.
Thus, we define the simplicial complex $K(L)$ to be the simplicial realization of $\K(L)$.
It follows that any $k$-simplex of $K(L)$ is given by a set of vertices of the form 
\begin{equation*}
    \{[(T_1,B,\lambda_1)], [(T_2,B,\lambda_2)], \ldots, [(T_{k+1},B,\lambda_{k+1})]\}
\end{equation*}
where $(T_i,B,\lambda_i) \leq (T_{i+1},B,\lambda_{i+1})$ \cite[\S 2.1]{mccullough-miller}.
Note that any simplex contains at most one labelled bipartite tree of any given type.
The group $\Aut(L)$ thus acts on $K(L)$ by simplicial automorphisms.
Observe that for every $\T$, any inner automorphism of $L$ is a product of vertex automorphisms carried by $\T$.
Thus, the action of $\Aut(L)$ on $K(L)$ factors through $\Out(L)$.

Finally, we note that $K(L)$ is not necessarily connected in general.
We let $K_0(L)$ denote the connected component containing the trivial $B_0$-labelled bipartite tree, where $B_0$ is our fixed preferred basis of $L$.
We say that two bases $B_1$ and $B_2$ are \emph{symmetrically equivalent} if each element of $B_2$ is a conjugate of an element of $B_1$.
Then a vertex $[(T,B,\lambda)]$ will lie in $K_0(L)$ if and only if $B$ is symmetrically equivalent to $B_0$, and the $\Out(L)$-stabilizer of the component $K_0(L)$ is $\SymOut(L)$ \cite[Prop~2.3]{mccullough-miller}.

\subsection{Relating the complexes}\label{subsec:complexes_quotient_map}

Now, fix a preferred basis $A_0 = \{y_1, \ldots, y_{n}\}$ of $F_{n}$, and a basis $B_0 = \{z_1, \ldots, z_{n}\}$ of $H_n$, and consider the projection $\pi:F_{n} \rightarrow H_n$ mapping $y_i$ to $z_i$.
Given a bipartite tree labelled by a basis $A$ of $F_n$, we can apply $\pi$ to each label to obtain a bipartite tree labelled by the set $\pi(A)$ (we warn the reader that $\pi(A)$ might not be a basis of $H_n$ for an arbitrary choice of $A$, but it will be a basis if $A$ is symmetrically equivalent to $A_0$).
Our next goal is to show that this yields well-defined quotient map $Q:K_0(F_n) \rightarrow K_0(H_n)$ which factors through an isomorphism $K_0(F_n)/\Ker(\P_n) \cong K_0(H_n)$.

\begin{lemma}
    The map $Q:K_0(F_n) \rightarrow K_0(H_n)$ is a well-defined equivariant map of simplicial complexes.
\end{lemma}
\begin{proof}
    Let $\widehat{\P}_n:\SymAut(F_n) \rightarrow \SymAut(H_n)$ be the lift of $\P_n$.
    Let $\widehat{\K}_0(F_n)$ be the set of isomorphism classes of labelled bipartite trees $(T, A, \lambda)$, where $A$ is symmetrically equivalent to $A_0$ (we are not quotienting by the relation $\sim$).
    Similarly, let $\widehat{\K}_0(H_n)$ be the set of isomorphism classes of bipartite trees labelled by a basis of $H_n$ which is symmetrically equivalent to $B_0$.
    
    First, we can check that there is a well-defined map 
    \begin{equation*}
        \widehat{Q}:\widehat{\K}_0(F_n) \rightarrow \widehat{\K}_0(H_n).
    \end{equation*}
    For $\T \in \widehat{K}_0(F_n)$, the labelled bipartite tree $\widehat{Q}(\T)$ is obtained by applying $\pi$ to each label of $\T$.
    To see that this is well-defined, it's enough to show that if $A$ is a basis of $F_n$ symmetrically equivalent to $A_0$, then $\pi(A)$ is a basis of $H_n$ symmetrically equivalent to $B_0$.
    Observe first that $A = f(A_0)$ for some symmetric automorphism $f$ of $F_n$.
    Indeed, since $A$ is a basis, there is some automorphism $f$ taking $A_0$ to $A$, and since $A$ is symmetrically equivalent to $A_0$, the automorphism $f$ must be symmetric.
    Then, 
    \begin{equation*}
        \pi(A) = \pi(f(A_0)) = \widehat{\P}_n(f)(\pi(A_0)) = \widehat{\P}_n(f)(B_0)
    \end{equation*}
    as desired.

    Next, observe that $\widehat{Q}$ descends to a well-defined map $Q:K_0(F_n) \rightarrow K_0(H_n)$.
    Indeed, if $f_{v,\theta} \in \Aut(F_n)$ is a vertex automorphism carried by $\T$, then $\widehat{\P}_n(f_{v,\theta})$ will be a vertex automorphism carried by $Q(\T)$, so if $\T \sim \T'$, then $\widehat{Q}(\T) \sim \widehat{Q}(\T')$.
    Moreover, since $\widehat{Q}$ preserves the type of a labelled bipartite tree, it follows that $\widehat{Q}$ preserves the partial ordering, and hence $Q$ is a map of simplicial complexes.

    Finally, since $\pi \circ f = \widehat{\P}_n(f) \circ \pi$ for any $f \in \SymAut(F_n)$, we see that $\widehat{Q}$ is equivariant with respect to the actions of $\SymAut(F_n)$ and $\SymAut(H_n)$, so it follows that $Q$ is equivariant with respect to the actions of $\SymOut(F_n)$ and $\SymOut(H_n)$.
\end{proof}

For the proof of the next lemma, we call some facts from \cite{mccullough-miller}.
Let $L$ denote $F_{n}$ or $H_n$.
Let $N(L)$ denote the set of vertices of $K_0(L)$ represented by trivial labelled bipartite trees (McCullough-Miller call these \emph{nuclear vertices}).
Given $v \in N(L)$, let $\st(v)$ denote the (closed) star of $v$ in $K_0(L)$.
Given a basis $B$ of $L$, let $\mathcal{T}_B$ be the poset of bipartite trees labelled by the (fixed) basis $B$ (McCullough-Miller call this the \emph{Whitehead poset}; its isomorphism class as a poset depends only on the size of $B$). We let $\v \mathcal{T}_B \v$ denote the simplicial realization of $\mathcal{T}_B$.
McCullough-Miller show the following \cite[\S 2.1]{mccullough-miller}:
\begin{enumerate}[label=(\roman*)]
    \item $\SymOut(L)$ acts transitively on $N(L)$,
    \item two trivial labelled bipartite trees represent the same vertex of $K_0(L)$ if and only if they differ by an inner automorphism,
    \item if $v \in N(L)$ is represented by the trivial $B$-labelled bipartite tree, then the natural map $\v \mathcal{T}_B \v \rightarrow \st(v)$ is a simplicial isomorphism,
    \item $K_0(L) = \bigcup_{v \in N(L)} \st(v)$.
\end{enumerate}

\begin{lemma}\label{lem:quotient_of_complexes}
    The map $Q:K_0(F_{n}) \rightarrow K_0(H_n)$ is surjective and factors through an isomorphism
    \begin{equation*}
        K_0(F_{n})/\Ker(\P_n) \cong K_0(H_n).
    \end{equation*}
\end{lemma}
\begin{proof}
    First, we claim that $Q$ restricts to an isomorphism $\st(v) \cong \st(Q(v))$ for each $v \in N(F_{n})$.
    Suppose $v$ is represented by an $A$-labelled trivial bipartite tree.
    Then $\pi(v)$ is represented by a $\pi(A)$-labelled trivial bipartite tree.
    The posets $\mathcal{T}_A$ and $\mathcal{T}_{\pi(A)}$ are naturally isomorphic, and we have a commuting square 
    \begin{equation*}
        \begin{tikzcd}
            {\v \mathcal{T}_A \v} & {\st(v)} \\
            {\v \mathcal{T}_{\pi(A)} \v} & {\st(\pi(v))}
            \arrow[from=1-1, to=1-2]
            \arrow["\cong", from=1-1, to=2-1]
            \arrow[from=1-2, to=2-2]
            \arrow[from=2-1, to=2-2]
        \end{tikzcd}
    \end{equation*}
    So, the claim follows from (iii).

    Next, we can prove that $Q$ is surjective.
    By (iv), and the fact that $Q$ restricts to an isomorphism on $\st(v)$ for $v \in N(F_n)$, it's enough to show that the restriction $Q_N:N(F_{n}) \rightarrow N(H_n)$ is surjective.
    Take any vertex $w \in N(H_n)$.
    Let $v_0 \in N(F_n)$ and $w_0 \in N(H_n)$ be the vertices labelled by $A_0$ and $B_0$ respectively, so $Q(v_0) = w_0$.
    By (iv), there is a sequence of adjacent vertices $w_0, u_0, w_1, u_1, \ldots, w_n=w$, where $w_i \in N(H_n)$ and $u_i \in \st(w_i) \cap \st(w_{i+1})$.
    We claim that we can lift this to a sequence of vertices $v_0, \hat{u}_0, v_1, \hat{u}_1, \ldots, v_n$ in $K_0(F_n)$ where $Q(v_i) = w_i$ and $Q(\hat{u}_i) = u_i$; from this claim, we get that $Q(v_n) = w$, and hence $Q_N$ is surjective.
    We start by lifting $w_0$ to $v_0$, and then proceeding inductively.
    If we can lift $w_i$ to $v_i$, then since $Q$ is an isomorphism on $\st(v_i)$, we can lift $u_i$ to some $\hat{u}_i$ adjacent to $v_i$.
    If we suppose $w_i$ is labelled by the basis $B_i$, then $B_i = \pi(A_i)$ for some basis $A_i$ labelling $v_i$.
    The vertex $u_i$ is then represented by a tree of the form $\T =(T, B_i, \lambda)$, and $\hat{u}_i$ is represented by a tree of the form $\widehat{\T} = ({T}, A_i, {\lambda})$.
    Suppose $w_{i+1}$ is labelled by the basis $B_{i+1}$. 
    Then since $u_i$ is adjacent to $w_{i+1}$, there must be some product of vertex automorphisms carried by $\T$ taking $B_i$ to $B_{i+1}$.
    But any vertex automorphism carried by $\T$ lifts to a vertex automorphism carried by $\widehat{\T}$.
    Thus we can find a lift $v_{i+1}$ of $w_{i+1}$ adjacent to $\hat{u}_i$.

    Next, we claim that the restriction of $Q$ to $Q_N:N(F_{n}) \rightarrow N(H_n)$ factors through a bijection 
    \begin{equation*}
        N(F_{n})/\Ker(\P_n) \leftrightarrow N(H_n).
    \end{equation*}
    First, we can show that $\Ker(\P_n)$ preserves each fiber of the map $Q_N$.
    We can see that for any $v = [(T,A,\lambda)] \in N(F_{n})$ and $\alpha = [f] \in \Ker(\P_n)$, the bases $A$ and $f(A)$ project to conjugate bases of $H_n$ (since $f$ descends to an inner automorphism of $H_n$), and hence $v$ and $\alpha(v)$ map to the same vertex of $K_0(H_n)$.
    Next, we can show that if $v,w \in N(F_{n})$ have the same image under $Q_N$, then they differ by an element of $\Ker(\P_n)$.
    Let $(T_v,A_v,\lambda_v)$ and $(T_w,A_w,\lambda_w)$ be trivial labelled bipartite trees representing $v$ and $w$ respectively.
    By (i), we can choose $\alpha \in \SymOut(F_{n})$ taking $v$ to $w$, and by (ii), we can choose a representative $f$ of $\alpha$ taking $A_v$ to $A_w$.
    Then $f$ descends to an automorphism $\overline{f}$ of $H_n$ taking $\pi(A_v)$ to $\pi(A_w)$.
    Since $Q(v) = Q(w)$, (ii) implies that $\pi(A_v)$ and $\pi(A_w)$ are conjugate bases, meaning that $\overline{f}$ is an inner automorphism, and thus $\alpha \in \Ker(\P_n)$.
    This proves the claim.

    Finally, we can prove that $Q$ factors through an isomorphism $K_0(F_n)/\Ker(\P_n) \cong K_0(H_n)$.
    Using the same argument as in the previous paragraph, we can see that for any vertex $v$ of $K_0(F_{n})$, $Q(v) = Q(\alpha(v))$ for any $\alpha \in \Ker(\P_n)$.
    So, it remains to show that if $v,w$ are two vertices of $K_0(F_{n})$ such that $Q(v) = Q(w)$, then $v$ and $w$ differ by an element of $\Ker(\P_n)$.
    Let $n_v,n_w \in N(F_{n})$ be vertices to which $v$ and $w$ are adjacent (these exist by (iv)).
    Since we have a bijection $N(F_{n})/\Ker(\P_n) \leftrightarrow N(H_n)$, we can choose $\alpha \in \Ker(\P_n)$ taking $n_v$ to $n_w$.
    Then $\alpha(v)$ is in $\st(n_w)$.
    Since $Q(\alpha(v)) = Q(v) = Q(w)$ and $Q$ is injective on $\st(n_w)$, it must be that $\alpha(v) = w$.    
\end{proof}

\subsection{Proof of Proposition \ref{prop:alg_kernel_computation}}\label{subsec:pf_alg_kernel_complutation}

Finally, we can prove Proposition \ref{prop:alg_kernel_computation}.
The main result of \cite{mccullough-miller} is that the complexes $K_0(F_n)$ and $K_0(H_n)$ are contractible.
By Lemma \ref{lem:quotient_of_complexes}, $\Ker(\P_n)$ acts on the contractible complex $K_0(F_n)$ and the quotient is contractible.
Moreover, we claim that $\Ker(\P_n)$ acts without rotations, meaning that if $\alpha \in \Ker(\P_n)$ preserves a simplex $s$ of $K_0(F_n)$, then $\alpha$ fixes each vertex of $s$. 
Indeed, as we remarked in Section \ref{subsec:the_complexes}, any simplex contains at most one labelled bipartite tree of any given type; since $\SymOut(F_n)$ preserves the type of a labelled bipartite tree, this implies that $\SymOut(F_n)$ acts without rotations, and hence $\Ker(\P_n)$ does too.
Thus, a result of Armstrong \cite{armstrong} says that $\Ker(\P_n)$ is generated by its stabilizers of vertices of $K_0(F_n)$ (see also \cite{putman-choices}).
This implies $\Ker(\P_n)$ is normally generated over $\SymOut(F_n)$ by the $\Ker(\P_n)$-stabilizers of one vertex from each $\SymOut(F_n)$-orbit.

Each $\SymOut(F_n)$-orbit contains a vertex $v$ represented by a labelled bipartite tree of the form $\T = (T, A_0, \lambda)$ where $A_0 = \{y_1, \ldots, y_n\}$ is our preferred basis.
McCullough-Miller \cite[Prop~5.1]{mccullough-miller} compute that the $\SymOut(F_n)$-stabilizer of such a vertex is 
\begin{equation*}
    \Gamma_v \coloneqq V_{\T} \rtimes ((\Z/2\Z)^n \rtimes S_\T)
\end{equation*}
where $V_{\T}$ is the group of vertex automorphisms carried by $\T$, $(\Z/2\Z)^n$ is the group generated by the elements $\rho_i$, and $S_\T$ is the group of permutations of $A_0$ that perserve $\T$ up to isomorphism.
Moreover, McCullough-Miller show that the group $V_{\T}$ is isomorphic to $\prod_{i=1}^n \Z^{k_i-1}$ where $k_i$ is the number of components of $T - \lambda(y_i)$.
Each $\Z^{k_i-1}$ factor corresponds to the group of vertex automorphisms for a fixed vertex $v = \lambda(y_i)$.
This is because a vertex automorphism $f_{v,\theta}$ amounts to a choice of power of $y_i$ for each component of $T - \lambda(y_i)$, and up to an inner automorphism we can always assume $f_{v,\theta}$ fixes one component.
The fact that we have a direct product of these $\Z^{k_i-1}$ factors is because any two vertex automorphisms based at different vertices commute up to inner automorphisms \cite[Lem~1.1]{mccullough-miller}.

The $\Ker(\P_n)$-stabilizer of a vertex is the intersection of $\Ker(\P_n)$ with the $\SymOut(F_n)$-stabilizer.
Thus, it remains to compute the intersection of $\Ker(\P_n)$ with $\Gamma_v$.
A direct computation shows that each $\rho_i$ is in $\Ker(\P_n)$, and hence the normal closure of $\{\rho_1, \ldots, \rho_n\}$ is contained in $\Ker(\P_n)$.
So, we want to show that $\Ker(\P_n) \cap \Gamma_v$ is contained in the normal closure of the set $\{\rho_1, \ldots, \rho_n\}$.
The images $\P_n(V_\T)$, $\P_n((\Z/2\Z)^n)$, and $\P_n(S_\T)$ pairwise intersect trivially, so it's enough to compute the intersection of each group with $\Ker(\P_n)$.
We see that the full group $(\Z/2\Z)^n = \la \rho_1, \ldots, \rho_n\ra$ lies in $\Ker(\P_n)$, and that the intersection $S_\T \cap \Ker(\P_n)$ is trivial.
Thus, it remains to compute $V_\T \cap \Ker(\P_n)$.
The restriction of $\P_n$ to $V_\T$ corresponds to the surjection
\begin{equation*}
    \prod_{i=1}^n \Z^{k_i-1} \rightarrow \prod_{i=1}^n (\Z/2\Z)^{k_i-1}.
\end{equation*}
So, $V_\T \cap \Ker(\P_n)$ is precisely $\prod_{i=1}^{n} (2\Z)^{k_i-1}$.
We therefore just need to show that for any $\alpha \in V_\T$, $\alpha^2$ lies in the normal closure of $\{\rho_1, \ldots, \rho_n\}$.
Suppose $\alpha$ represented by a vertex automorphism $f_{v,\theta}$ with $v = \lambda(y_i)$.
Then we can compute that
\begin{equation*}
    \rho_i f_{v,\theta} \rho_i^{-1} = f_{v,\theta}^{-1}.
\end{equation*}
This implies $f_{v,\theta}^2 = [f_{v,\theta},\rho_i]$, and hence $\alpha^2$ lies in the normal closure of the set $\{\rho_1, \ldots, \rho_n\}$ as desired.

\section{Computing the hyperelliptic lifting kernel}\label{sec:pf_thm_C}
In this section, we prove Theorem \ref{mainthm:hyperelliptic_kernel} and Corollary \ref{maincor:explicit_kernel_computation} by finding a normal generator of the lifting map
\begin{equation*}
    \L_n:\Mod(S^3,C_n) \rightarrow \SMod_n(M_n)/\la \tau \ra.
\end{equation*}
As before, we let $\rho_i \in \Mod(S^3, C_n) \cong \SymOut(F_n)$ be the element which inverts the $i$th generator of $\pi_1(S^3 - C_n) \cong F_n$.
Our goal is to show that $\Ker(\L_n) = \la \la \rho \ra \ra$, where $\rho = \prod_{i=1}^n \rho_i$ and $\la \la \rho \ra \ra$ is the $\Mod(S^3,C_n)$-normal closure of $\{\rho\}$.

Propositions \ref{prop:lifting_kernel_is_alg_kernel} and \ref{prop:alg_kernel_computation} tell us that $\Ker(\L_n)$ is contained in the normal closure of the set $\{\rho_1, \ldots, \rho_n\}$.
This means that any element of $\Ker(\L_n)$ is a product of conjugates of the elements $\rho_i$.
To complete the proof of Theorem \ref{mainthm:hyperelliptic_kernel}, we use the presentation of $\Mod(S^3,C_n) \cong \SymOut(F_n)$ in Section \ref{subsec:liftable_mcg_unlink} to conclude that any element of $\Ker(\L_n)$ is a product of conjugates of $\rho$.
Then, we will deduce Corollary \ref{maincor:explicit_kernel_computation} again using this presentation.

In Section \ref{subsec:kernel_lower_bound}, we show that $\la \la \rho \ra \ra \subseteq \Ker(\L_n)$.
In Section \ref{subsec:two_factor_case}, we examine the case $n=2$, and prove Theorem \ref{mainthm:hyperelliptic_kernel} and Corollary \ref{maincor:explicit_kernel_computation} directly for this case.
In Section \ref{subsec:semipalindromes}, we prove an auxiliary lemma which says that certain special words in $\Mod(S^3,C_n)$, called \emph{semipalindromes}, lie in $\la \la \rho \ra \ra$.
In Section \ref{subsec:pf_thm_C}, we complete the proof of Theorem \ref{mainthm:hyperelliptic_kernel} by showing that $\Ker(\L_n) \subseteq \la \la \rho \ra \ra$, and in Section \ref{subsec:pf_cor_D}, we complete the proof of Corollary \ref{maincor:explicit_kernel_computation} in the case $n=3$.

\subsection{Lower bound on the kernel}\label{subsec:kernel_lower_bound}

Our first step is to verify that $\la \la \rho \ra \ra \subseteq \Ker(\L_n)$.
In fact, we have the following precise statement.
Recall from our model of $p_n$ in Section \ref{subsec:model_of_cover} that $M_n$ can be built by gluing two copies $W_1$ and $W_2$ of $S^3 - (\text{$n$ open $3$-balls})$.
For $1 \leq i \leq n$, we let $\Sigma_i \subseteq M_n$ be the sphere coming from the $i$th boundary component of $W_1$ and $W_2$.

\begin{lemma}\label{lem:kernel_lower_bound}
    For $1 \leq i \leq n$, the element $\rho_i$ maps to the sphere twist $T_{\Sigma_i}$.
    Moreover, 
    \begin{equation*}
        \Ker(\L_n) \cap \la \rho_1, \ldots, \rho_n \ra = \la \rho \ra.
    \end{equation*}
\end{lemma}

In other words, $\L_n$ maps $\la \rho_1, \ldots, \rho_n\ra$ onto $\Twist(M_n)$, and the map $\la \rho_1, \ldots, \rho_n \ra \rightarrow \Twist(M_n)$ corresponds to the quotient map 
\begin{equation*}
    (\Z/2\Z)^n \rightarrow (\Z/2\Z)^n/\la e_1 + \cdots + e_n \ra \cong (\Z/2\Z)^{n-1},
\end{equation*}
where $e_1, \ldots, e_n$ is the natural basis of $(\Z/2\Z)^n$.

\begin{proof}[Proof of Lemma \ref{lem:kernel_lower_bound}]
    First, we can show that $\rho_i$ lifts to $T_{\Sigma_i}$.
    To do so, we'll construct an equivariant representative of $T_{\Sigma_i}$ that descends to a representative of $\rho_i$.
    Choose a regular neighborhood $U \cong S^2 \times [0,1]$ of $\Sigma_i$ where $\Sigma_i$ is the slice $S^2 \times \{1/2\}$.
    We can choose this neighborhood so that $S^2 \times [0,1/2] \subseteq W_1$ and $S^2 \times [1/2,1] \subseteq W_2$, and so that $\tau$ takes $S^2 \times \{t\}$ to $S^2 \times \{1 - t\}$ (so $U$ is $\tau$-invariant).
    Recall that we view $S^3$ as $\R^3 \cup \{\infty\}$ and assume $\Sigma_i$ is centered on the $x$-axis in both $W_1$ and $W_2$.
    Let $f$ be the homeomorphism of $M_n$ that acts as follows:
    \begin{itemize}
        \item On $S^2 \times [0,1/2] \subseteq W_1$, $f$ applies a $2\pi t$ rotation about the $x$-axis on the slice $S^2 \times \{t\}$.
        \item On $S^2 \times [1/2,1] \subseteq W_2$, $f$ applies a $-2\pi t$ rotation about the $x$-axis on the slice $S^2 \times \{t\}$.
        \item Outside of $U$, $f$ acts as the identity map.
    \end{itemize}
    Since $\tau$ swaps $S^2 \times \{t\}$ with $S^2 \times \{1-t\}$ and a $-2\pi(1-t)$ rotation is a $2\pi t$ rotation, we see that $f$ commutes with $\tau$, i.e.\ $f$ is symmetric.
    Moreover, if we let $K_i$ denote the $i$th component of $C_n$ (i.e.\ the component of $C_n$ coming from $\Sigma_i$), then $f$ descends to a homeomorphism $\overline{f}$ which is supported in a neighborhood of $p_n(\Sigma_i)$ and reverses the orientation of $K_i$.
    It follows that $\overline{f}$ must represent the element $\rho_i$.

    Recall that $\Twist(M_n) \cong (\Z/2\Z)^{n-1}$ and is generated by the sphere twists about the spheres $\Sigma_1, \ldots, \Sigma_{n-1}$.
    Moreover, the sphere twists $T_{\Sigma_i}$ satisfy the relation 
    \begin{equation*}
        \prod_{i=1}^n T_{\Sigma_i} = 1
    \end{equation*}
    (see e.g.\ \cite[\S 2]{hatcher-wahl}).
    Thus the kernel of $\L_n$ restricted to $\la \rho_1, \ldots, \rho_n \ra$ is precisely $\rho$.
\end{proof}

\subsection{The case of two factors}\label{subsec:two_factor_case}
Next, we can directly prove Corollary \ref{maincor:explicit_kernel_computation} in the case $n=2$.
For $n=2$, the elements $\alpha_{i,j}$ in our presentation of $\Mod(S^3,C_n)$ (see Section \ref{subsec:liftable_mcg_unlink}) are trivial.
Thus, $\Mod(S^3, C_2)$ is generated by order $2$ elements $\rho_1$, $\rho_2$, and $\sigma$, where $\sigma$ swaps the components of $C_2$.
Further, we see that 
\begin{equation*}
    \Mod(S^3, C_2) \cong (\Z/2\Z \times \Z/2\Z) \rtimes \Z/2\Z.
\end{equation*}
On the other hand, we have that $M_2 = S^1 \times S^2$, and Gluck \cite{gluck} proved that 
\begin{equation*}
    \Mod(S^1 \times S^2) \cong \Z/2\Z \times \Z/2\Z
\end{equation*}
where the first factor is a sphere twist $T$ about the core sphere, and the second factor is represented by $\tau$ (this is a special case of the results of \cite{laudenbach-1,laudenbach-2} and \cite{brendle-broaddus-putman}).

From Lemma \ref{lem:kernel_lower_bound}, we see that $\rho_1$ and $\rho_2$ both lift to $T$.
We claim that $\sigma$ lifts to a homeomorphism which is isotopic to $\tau$.
Indeed, let $\theta \in \Homeo^+(S^1 \times S^2)$ be a $\pi$ rotation in the $S^1$ factor, and let $\tau' = \theta \circ \tau$.
Then $\tau'$ is symmetric since $\theta$ commutes with $\tau$, and it descends to a homeomorphism of $(S^3, C_2)$ which permutes the components of $C_2$, which must represent $\sigma$.

This means that every element of $\Mod(S^1 \times S^2)$ has symmetric representative, i.e.\ $\SMod(S^1 \times S^2) = \Mod(S^1 \times S^2)$ and $\SMod(S^1 \times S^2)/ \la \tau \ra \cong \Z/2\Z$.
We see that $\Ker(\L_2$) is precisely $\la \rho, \sigma \ra \cong \Z/2\Z \times \Z/2\Z$.

\subsection{Semipalindromes}\label{subsec:semipalindromes}
Now, we can begin to prove Theorem \ref{mainthm:hyperelliptic_kernel} and Corollary \ref{maincor:explicit_kernel_computation} for the case $n \geq 3$.
We start with a group-theoretic lemma about the generators of $\Mod(S^3, C_n)$.
As in Section \ref{subsec:liftable_mcg_unlink}, for $1 \leq i,j \leq n$ we let $\alpha_{i,j} \in \Mod(S^3, C_{n+1})$ be the element that maps $y_i$ to $y_jy_iy_j^{-1}$ and fixes each other generator.

We define a \emph{semipalindrome} to be a word on the letters $\alpha_{i,j}^{\pm 1}$ that can be obtained by the following inductive rules:
\begin{enumerate}[label=(\roman*)]
    \item $\id$ and $\alpha_{i,j}^{\pm 2}$ are semipalindromes
    \item if $w$ is a semipalindrome, then so are $\alpha w \alpha$ and $\alpha w\alpha^{-1}$ for any $\alpha = \alpha_{i,j}^{\pm 1}$.
\end{enumerate}
Because $\rho_k \alpha_{i,j} \rho_k^{-1} = \alpha_{i,j}^{\pm 1}$, we see that for any $1 \leq k \leq n$ and $v \in \PMod(S^3, C_n)$, $v\rho_kv^{-1} = w\rho_k$ for some semipalindrome $w$.
Moreover, semipalindromes have the following key property.

\begin{lemma}\label{lem:semipalindrome_in_normal_closure}
    Any semipalindrome lies in $\la \la \rho \ra \ra$.
\end{lemma}
\begin{proof}
    We proceed by induction.
    Observe first that since $\rho\alpha_{i,j}\rho^{-1} = \alpha_{i,j}^{-1}$, we have that
    \begin{equation*}
        \alpha_{i,j}^{\pm 2} = [\alpha_{i,j}^{\pm 1}, \rho] \in \la \la \rho \ra \ra.
    \end{equation*}
    Now, let $w$ be a semipalindrome that lies in $\la \la \rho \ra \ra$.
    Let $\alpha$ denote some $\alpha_{i,j}$ or $\alpha_{i,j}^{-1}$.
    Then $\alpha w \alpha^{-1}$ automatically lies in $\la \la \rho \ra \ra$, so we just have to show that $\alpha w \alpha$ lies in $\la \la \rho \ra \ra$.

    On the one hand, since $w \in \la \la \rho \ra \ra$, we can write
    \begin{equation*}
        w = \prod_{k=1}^m v_k\rho v_k^{-1}
    \end{equation*}
    for some $v_k \in \Mod(S^3, C_n)$.
    In fact, using the presentation in Section \ref{subsec:liftable_mcg_unlink}, we may assume that each $v_k$ lies in $\PMod(S^3, C_n)$, i.e.\ $v_k$ is a product of $\alpha_{i,j}^{\pm 1}$'s.
    From the relation $\rho\alpha_{i,j}\rho^{-1} = \alpha_{i,j}^{-1}$, we can write
    \begin{equation*}
        v_k \rho v_k^{-1} = w_k \rho
    \end{equation*}
    for some semipalindrome $w_k$.
    Moreover, we can write
    \begin{equation*}
        w = \prod_{k=1}^m v_k\rho v_k^{-1}
        = \prod_{k=1}^m w_k\rho 
        = \left( \prod_{k=1}^m w_k' \right) \rho^\ell
    \end{equation*}
    where $w_k'$ is the semipalindrome $w_k$ or $\rho w_k \rho^{-1}$, and $\ell$ is either $0$ or $1$ depending on the parity of $m$.
    Note that $w_k$ lies in $\la \la \rho \ra \ra$, since
    \begin{equation*}
        w_k = w_k\rho^2 = v_k\rho v_k^{-1} \rho \in \la \la \rho \ra \ra.
    \end{equation*}
    It follows that $w_k'$ also lies in $\la \la \rho \ra \ra$, and hence so does $\prod_{k=1}^m w_k'$
    
    Now, if $\ell = 1$, we have 
    \begin{equation*}
        \alpha w \alpha
        = \alpha \left( \prod_{k=1}^m w_k' \right) \rho \alpha 
        = \alpha \left( \prod_{k=1}^m w_k' \right) \alpha^{-1} \rho \in \la \la \rho \ra \ra.
    \end{equation*}
    Otherwise, if $\ell = 0$, we can write 
    \begin{equation*}
        \alpha w \alpha
        = \alpha \left( \prod_{k=1}^m w_k' \right) \alpha 
        = \alpha \left( \prod_{k=1}^m w_k'\alpha^{-1}\alpha \right) \alpha 
        = \left( \prod_{k=1}^m \alpha w_k' \alpha^{-1} \right) \alpha^2 \in \la \la \rho \ra \ra.
    \end{equation*}
\end{proof}

\subsection{Proof of Theorem \ref{mainthm:hyperelliptic_kernel}}\label{subsec:pf_thm_C}

Finally, we can prove Theorem \ref{mainthm:hyperelliptic_kernel}.
Our goal is to show that $\Ker(\L_n) = \la \la \rho \ra \ra$.
Lemma \ref{lem:kernel_lower_bound} tells us that $\la\la \rho \ra \ra \subseteq \Ker(\L_n)$, so it remains to show that $\Ker(\L_n) \subseteq \la \la \rho \ra \ra$.

First, we claim that for any $\beta \in \Ker(\L_n)$, we can write
\begin{equation*}
    \beta = \left(\prod_{k=1}^m w_k \right) \left( \prod_{i=1}^n \rho_i^{\ell_i} \right)
\end{equation*}
for some semipalindromes $w_k$ and some $\ell_i \geq 1$.
To prove this, observe first that by Propositions \ref{prop:lifting_kernel_is_alg_kernel} and \ref{prop:alg_kernel_computation}, we know that $\Ker(\L_n) \subseteq \la \la \rho_1, \ldots, \rho_n \ra \ra$.
So, we can write
\begin{equation*}
    \beta = \prod_{k=1}^m v_k \rho_{i_k} v_k^{-1}
\end{equation*}
for some $v_k \in \Mod(S^3, C_n)$.
From the relations $\rho_k \alpha_{i,j} \rho_k^{-1} = \alpha_{i,j}^{\pm 1}$ and $\sigma \alpha_{i,j} \sigma^{-1} = \alpha_{\sigma(i), \sigma(j)}$ for $\sigma \in S_n$, we can assume each $v_k$ in fact lies in $\PMod(S^3, C_n)$ (i.e.\ each $v_k$ is a word in the elements $\alpha_{i,j}^{\pm 1}$).
The claim then follows by applying the relation $\rho_k \alpha_{i,j} \rho_k^{-1} = \alpha_{i,j}^{\pm 1}$.

Now, by Lemmas \ref{lem:kernel_lower_bound} and \ref{lem:semipalindrome_in_normal_closure}, the map $\L_n$ takes each $w_k$ to $\id$.
Since $\L_n(\beta) = \id$, it follows from Lemma \ref{lem:kernel_lower_bound} that $\ell_1 = \cdots = \ell_n$, i.e.\ that 
\begin{equation*}
    \beta = \left( \prod_{k=1}^m w_k \right) \rho^\ell
\end{equation*}
for some $\ell \geq 0$.
Finally, we conclude from Lemma \ref{lem:semipalindrome_in_normal_closure} that $\beta \in \la \la \rho \ra \ra$.

\subsection{Proof of Corollary \ref{maincor:explicit_kernel_computation}}\label{subsec:pf_cor_D}

Our final task is to prove Corollary \ref{maincor:explicit_kernel_computation} for the case $n=3$.
We want to show that $\Ker(\L_3) \cong F_\infty \rtimes \Z/2\Z$.

Recall from Section \ref{subsec:liftable_mcg_unlink} that
\begin{equation*}
    \Mod(S^3,C_3) \cong \PMod(S^3,C_3) \rtimes ((\Z/2\Z)^3 \rtimes S_3).
\end{equation*}
By Theorem \ref{mainthm:hyperelliptic_kernel}, $\Ker(\L_3) = \la \la \rho \ra \ra$.
A direct computation shows that this subgroup is contained in $\PMod(S^3,C_3) \rtimes \la \rho \ra$.

From Section \ref{subsec:liftable_mcg_unlink}, the group $\PMod(S^3,C_3)$ is generated by the elements $\alpha_{i,j}$ for $1 \leq i, j \leq 3$ (where $\alpha_{i,i}$ is trivial) subject to the relations
\begin{align*}
    [\alpha_{i,j},\alpha_{k,\ell}] &= 1 \\
    [\alpha_{i,k}, \alpha_{j,k}] &= 1 \\
    \alpha_{i,j}\alpha_{j,k}\alpha_{i,k} &= \alpha_{i,k}\alpha_{j,k}\alpha_{i,j}, \\
    \prod_{i=1}^3 \alpha_{i,j} &= 1,
\end{align*}
where $i$, $j$, $k$, and $\ell$ are distinct indices.
The first relation is vacuous since there is no possible choice of four distinct indices.
The second and third relations are redundant by the fourth relation, since $\alpha_{i,k} = \alpha_{j,k}^{-1}$.
Thus, $\PMod(S^3,C_n)$ is freely generated by the elements $\alpha_{1,2}$, $\alpha_{2,3}$, and $\alpha_{3,1}$.
It follows that $\PMod(S^3,C_3) \rtimes \la \rho \ra$ is generated by the elements $\alpha_{1,2}$, $\alpha_{2,3}$, $\alpha_{3,1}$, and $\rho$, subject to the relations $\rho^2 = 1$ and $\rho \alpha_{i,j} \rho^{-1} = \alpha_{i,j}^{-1}$.
That is,
\begin{equation*}
    \PMod(S^3,C_3) \rtimes \la \rho \ra \cong F_3 \rtimes \Z/2\Z.
\end{equation*}

If we quotient $\PMod(S^3,C_3) \rtimes \la \rho \ra$ by $\la \la \rho \ra \ra$, this amounts to adding the relation $\rho = 1$.
The resulting group will be generated by the elements $\overline{\alpha}_{1,2}$, $\overline{\alpha}_{2,3}$ and $\overline{\alpha}_{3,1}$ subject to the relations $\overline{\alpha}_{i,j}^2 = 1$, which is isomorphic to the group $H_3 = (\Z/2\Z) * (\Z/2\Z) * (\Z/2\Z)$.
In other words, $\la \la \rho \ra \ra$ is isomorphic to the kernel of the map 
\begin{equation*}
    \omega: F_3 \rtimes \Z/2\Z \rightarrow H_3
\end{equation*}
which maps the $\Z/2\Z$ factor to $1$ and reduces each generator of $F_3$ mod $2$.
If we let $\pi:F_3 \rightarrow H_3$ be the natural projection, then $\Ker(\omega) \cong \Ker(\pi) \rtimes \Z/2\Z$.
Corollary \ref{maincor:explicit_kernel_computation} then follows from the fact that $\Ker(\pi)$ is an infinite rank free group.
 
\printbibliography

\end{document}